\documentclass[reqno]{amsart}

\usepackage{amssymb,amsmath}
\usepackage{euscript}
\usepackage{graphicx}
\usepackage{tikz}
\usepackage{tikz-cd}
\usepackage{hyperref}
\usepackage{mathrsfs}
\usepackage{dsfont}
\usepackage{mathtools}
\usepackage{array}
\usepackage{tabu}

\newtheorem{thm}{Theorem}[section]
\newtheorem{cor}[thm]{Corollary}
\newtheorem{lem}[thm]{Lemma}
\newtheorem{prop}[thm]{Proposition}
\newtheorem{defin}[thm]{Definition}
\newtheorem{def-lem}[thm]{Definition-Lemma}

\theoremstyle{remark}

\newtheorem{rem}[thm]{Remark}

\newtheorem{assu}{Assumption}

\numberwithin{equation}{section}

\newcommand{\bbC}{\mathbb{C}}

\newcommand{\bbN}{\mathbb{N}}

\newcommand{\bbR}{\mathbb{R}}

\newcommand{\bbZ}{\mathbb{Z}}

\newcommand{\scrL}{\mathscr{L}}
\newcommand{\scrM}{\mathscr{M}}

\newcommand{\scrU}{\mathscr{U}}
\newcommand{\scrV}{\mathscr{V}}
\newcommand{\scrW}{\mathscr{W}}

\newcommand{\calC}{\mathcal{C}}

\newcommand{\calL}{\mathcal{L}}

\newcommand{\frakS}{\mathfrak{S}}

\DeclareMathOperator{\crit}{Crit}
\DeclareMathOperator{\grad}{grad}

\newcommand{\cs}{\mathrm{cs}}
\newcommand{\edgef}{\mathrm{Edge_{\mathrm{int}}}}
\newcommand{\edgesi}{\mathrm{Edge}_\infty}

\newcommand{\pre}{\mathrm{pre}}

\newcommand{\vertex}{\mathrm{Vert}}

\newcommand{\abs}[1]{\lvert#1\rvert}

\title{Stable Morse flow trees}
\author{Kenneth Blakey}
\address{Department of Mathematics, MIT, 182 Memorial Drive, Cambridge, MA 02139, U.S.A.} 
\email{kblakey@mit.edu}

\begin{document}

\begin{abstract}
We prove a Floer-Gromov compactness type result for (stable) Morse flow trees of Legendrians in 1-jet spaces with simple front singularities satisfying Ekholm's preliminary transversality condition.
\end{abstract}

\maketitle
\tableofcontents

\section{Introduction}
\subsection{Context} 
Let $(M,g)$ be a closed Riemannian $n$-manifold, $J^0M\equiv M\times\bbR$ its 0-jet space, and $J^1M\equiv T^*M\times\bbR$ its 1-jet space. We endow $J^1M$ with its standard contact structure given as the kernel of the contact 1-form $dz-\alpha$, where $z$ is the $\bbR$-coordinate and $\alpha$ is the (pullback of the) tautological 1-form on $T^*M$. Recall, an $n$-dimensional submanifold of $J^1M$ is called \emph{Legendrian} if it is an integral submanifold of the ambient contact structure. Also, we endow $T^*M$ with its standard symplectic structure given as $-d\alpha$. Recall, an $n$-dimensional submanifold of $T^*M$ is called \emph{Lagrangian} if the restriction of the ambient symplectic form vanishes.

Legendrian contact homology is an isotopy invariant of Legendrians in the framework of symplectic field theory \cite{EGH00}; Ekholm's Morse flow trees were introduced in \cite{Ekh07} in order to compute the Legendrian contact homology of Legendrians in 1-jet spaces. In particular, the rigorous construction of Legendrian contact homology in 1-jet spaces was performed in \cite{EES07}, where the construction involves counting rigid pseudo-holomorphic disks in $T^*M$ with boundary on the projection of $L$ and punctures asymptotic to the projection of Reeb chords on $L$ -- finding pseudo-holomorphic disks involves solving a non-linear first order PDE, thus explicit computations are quite difficult. \cite[Theorem 1.1]{Ekh07} reduces this infinite-dimensional problem to a finite-dimensional problem in Morse theory by building a 1-1 correspondence between these rigid pseudo-holomorphic disks with at most 1 positive puncture and rigid Morse flow trees with exactly 1 positive puncture; thus, the differential of Legendrian contact homology in 1-jet spaces can be computed in terms of Morse theory.\footnote{At least, for Legendrians with simple front singularities, cf. Definition \ref{defin:simplefrontsingularities}.} The relation between Morse theory and pseudo-holomorphic curves can already be found in the Floer-theoretic case, cf. \cite{Flo89c,FO97}, and the power of Morse flow trees can be seen in various computations, cf. \cite{DR11,DRG23,EENS13,EN20}.

\subsection{Main results} 
Let $L\subset J^1M$ be a Legendrian with simple front singularities; moreover, we will assume $L$ satisfies the preliminary transversality condition detailed in \cite[Pages 1101-1103]{Ekh07}. We may locally describe $L$ as the multi 1-jet lift of smooth functions locally defined on $M$. A Morse flow tree of $L$ is essentially a finite tree decorated with the following additional data: for each edge, we assign a gradient flow of a difference of functions locally describing $L$, cf. Definition \ref{defin:morseflowtree}. In the present article, we generalize Morse flow trees to (pre)stable Morse flow trees essentially by allowing broken gradient flows, cf. Definitions \ref{defin:prestablemorseflowtree} \& \ref{defin:stablemorseflowtree}. We define a notion of (pre-)Floer-Gromov type convergence for (pre)stable Morse flow trees, cf. Definitions \ref{defin:prefloergromov} \& \ref{defin:floergromov}. Let $\frakS(L)$ be the moduli space of stable Morse flow trees of $L$. Let $N\in\bbZ_{>0}$ and consider the moduli space $\frakS(L,N)\subset\frakS(L)$ of $T$ such that the domain of $T$ has at most $N$ edges; we may define the moduli space 
    \begin{equation}
    \overline{\frakS}(L,N)\equiv\big\{T\in\frakS(L):\exists\{T_\nu\}\subset\frakS(L,N)\;\;\textrm{Floer-Gromov converging to $T$}\big\}.
    \end{equation}
Our main result is the following.

\begin{thm}\label{thm:main}
There exists a topology on $\overline{\frakS}(L,N)$ which is compact; moreover, this topology is unique with respect to the property that a sequence topologically converges if and only if that sequence Floer-Gromov converges.
\end{thm}

Theorem \ref{thm:main} is expected given Ekholm's correspondence and the corresponding statement of Floer-Gromov type convergence for pseudo-holomorphic disks.

\begin{rem}
As mentioned above, Ekholm's correspondence is for rigid pseudo-holomorphic disks with at most 1 positive puncture and rigid Morse flow trees with exactly 1 positive puncture. This correspondence fails, in general, for pseudo-holomorphic disks with many positive punctures; this is because one must study higher-dimensional moduli spaces of possibly multiply-covered Morse flow trees. (Recall, any Morse flow tree with exactly 1 positive puncture is not multiply-covered.) The present article allows possibly multiply-covered (pre)stable Morse flow trees, and we will not discuss their connection to pseudo-holomorphic disks.
\end{rem}

The main practical application of Theorem \ref{thm:main} is the following. Consider a sequence $\{T_\nu\}$ of Morse flow trees such that each $T_\nu$ has at most $P\in\bbZ_{>0}$ positive punctures (cf. Definition \ref{defin:puncture}) and at most formal dimension $D\in\bbZ_{>0}$ (cf. Definition \ref{defin:fd}). By \cite[Lemma 3.12]{Ekh07}, there exists a constant $k(M,L,P,D)\in\bbZ_{>0}$, depending on $M$, $L$, $P$, and $D$, such that the domain of each $T_\nu$ has at most $k(M,L,P,D)$ edges, i.e., the number of edges of a Morse flow tree can be bounded in terms of geometric data. We may now think of $\{T_\nu\}$ as a sequence in $\overline{\frakS}\big(L,k(M,L,P,D)\big)$ and construct a Floer-Gromov converging subsequence with limit $T$. The Floer-Gromov limit $T$ can exhibit two kinds of phenomena.
\begin{enumerate}
\item First, the gradient flows of a local function difference can break, in the standard way, at critical points of that local function difference.
\item Second, since the gradient flow assigned to an edge does not have to be defined on its maximal domain of definition, edges are allowed to contract to 0 length, i.e., a sequence may bubble off ``ghost edges''.
\end{enumerate}
In particular, the two aforementioned phenomena allow a sequence of stable Morse flow trees to develop nodal vertices, cf. Definitions \ref{defin:node1} \& \ref{defin:node2}; these are vertices at which a stable Morse flow tree can be decomposed into the gluing of two stable Morse flow trees, cf. Remarks \ref{rem:node1} \& \ref{rem:node2} for examples.

\subsection*{Acknowledgments}
The work for the present article was completed while the author was participating in the 2021 DIMACS REU. The author would like to thank Christopher Woodward for suggesting and mentoring this project during that time. The author would also like to thank Soham Chanda, Tyler Lane, Eric Kilgore, and Benoit Pausader for conversations. The author acknowledges support from the Rutgers University Department of Mathematics and NSF grant DMS-1711070 while completing the work for this project.

\section{Preliminaries on Morse flow trees}
\subsection{Metric ribbon trees} We begin by quickly reviewing the geometric realization of Stasheff's associahedron \cite{Sta63} given via stable metric ribbon trees -- this is originally due to Boardman-Vogt \cite{BV73}. We follow the presentation in \cite{FO97,MW10}.

Let $k\in\bbZ$ with $k\geq2$. Recall, the $k$-th associahedron $K_k$ is a CW-complex of dimension $k-2$ whose vertices correspond to the possible ways of inserting parentheses into the formal expression 
    \begin{equation}
    x_1\cdots x_k.
    \end{equation}
Each facet of $K_k$ is naturally the image of an embedding 
    \begin{equation}
    K_i\times K_j\to K_k.
    \end{equation}
    
\begin{defin}
A \emph{ribbon tree} consists the following data. 
\begin{enumerate}
\item A finite tree 
    \begin{equation}
    T\equiv\big(\vertex(T),\edgef(T)\big),
    \end{equation}
where $\vertex(T)$ is a set of vertices and $\edgef(T)$ a set of (interior) edges each incident to 2 vertices, together with a set $\edgesi(T)$ of semi-infinite edges each incident to 1 vertex.
\item For each $v\in\vertex(T)$, a cyclic ordering of the set 
    \begin{equation}
    \{e\in\edgef(T):v\in e\}.
    \end{equation}
\item A distinguished vertex $e_0\in\edgesi(T)$ called the \emph{root} (the other semi-infinite edges are referred to as \emph{leaves}).
\end{enumerate}
\end{defin}

A \emph{metric ribbon tree} is a ribbon tree $T$ together with a function 
    \begin{equation}
    d:\edgef(T)\to\bbR_{>0}
    \end{equation}
referred to as the \emph{metric}. We say a (metric) ribbon tree $T$ is \emph{stable} if, for each vertex $v\in\vertex(T)$, the sum of the valency of $v$ and
    \begin{equation}
    \big\lvert\{e\in\edgesi(T):v\in e\}\big\rvert
    \end{equation}
is at least 3. Given a stable ribbon tree $T$ with $k-1$ leaves, we denote by $G_{k,T}$ the set of all metrics on $T$; the space of all stable metric ribbon trees with $k-1$ leaves is denoted 
    \begin{equation}\label{eq:celldecomp}
    G_k\equiv\bigcup_{T}G_{k,T}
    \end{equation}
and has the following topology. Let $\{d_\nu\}\subset G_{k,T}$ be a sequence such that 
    \begin{equation}
    \lim_{\nu\to\infty} d_\nu(e)=d_\infty(e)\in\bbR_{\geq0}
    \end{equation}
for each $e\in\edgef(T)$. We define a stable ribbon tree $T_\infty$ by collapsing all edges $e\in\edgef(T)$ such that $d_\infty(e)=0$; we endow $T_\infty$ with the metric $\overline{d}_\infty$ given by the obvious restriction of $d_\infty$. We define the limit of $d_\nu$ to be $\overline{d}_\infty$. Now, \eqref{eq:celldecomp} provides a cell decomposition of $G_k$; in fact, $G_k$ is homeomorphic to $\bbR^{k-3}$.

Clearly, the closure of $G_{k,T}$ inside of $G_k$, in the topology just discussed, is given by $\cup_{T'\leq T}G_{k,T'}$, where $T'\leq T$ means $T'$ is obtained from $T$ by collapsing interior edges. Each cell $G_{k,T}$ may be compactified by allowing interior edges to have infinite length; we denote by $\overline{G}_k$ the induced compactification of $G_k$. Finally, Boardman-Vogt show that $\overline{G}_k$ is homeomorphic to $K_k$, and this homeomorphism is compatible with the various inclusions of facets.

\subsection{Local function differences}
We now return to Legendrians in 1-jet spaces following \cite[Section 2.2]{Ekh07}. As before, let $(M,g)$ be a closed Riemannian $n$-manifold and $L\subset J^1M$ a closed Legendrian in the 1-jet space of $M$. Via a generic perturbation of $L$, we may assume $L$ is \emph{chord generic}, i.e., $L$ only has finitely many Reeb chords. Consider the Lagrangian projection,
    \begin{equation}
    \Pi_\bbC:J^1M\to T^*M,
    \end{equation}
which forgets the $\bbR$-coordinate. Observe, $\Pi_\bbC(L)$ is an immersed Lagrangian. Since the Reeb field of the standard contact structure on $J^1M$ is $\partial_z$, it follows that Reeb chords on $L$ are in 1-1 correspondence with self-transverse double points of $\Pi_\bbC(L)$. We will assume $L$ satisfies the following.

\begin{defin}\label{defin:simplefrontsingularities}
We say that $L$ has \emph{simple front singularities} if the following conditions hold.
\begin{itemize}
\item The \emph{base projection} $\Pi:J^1M\to M$, when restricted to $L$, is an immersion outside of a closed codimension 1 submanifold $\Sigma\subset L$.
\item We denote by $\Pi_F:J^1M\to J^0M$ the \emph{front projection}. For any $\widetilde{p}\in\Sigma$, we have that $\Pi_F$, when restricted to $L$, has a standard \emph{cusp-edge singularity} at $\widetilde{p}$. Explicitly, this means there exists coordinates $y\equiv(y_1,\ldots,y_n)$ on $L$ near $\widetilde{p}$ and coordinates $x\equiv(x_1,\ldots,x_n)$ on $M$ near $\Pi(\widetilde{p})$ such that, if $z$ also denotes the $\bbR$-coordinate on $J^0M$, then 
    \begin{equation}
    \Pi_F(y)=\big(x(y),z(y)\big), 
    \end{equation}
where 
    \begin{align}\label{eqn:cuspedgecoordinates}
    x_1(y)&\equiv\dfrac{1}{2}y_1^2, \\
    x_i(y)&\equiv y_i,\;\;i\geq2, \nonumber \\
    z(y)&\equiv\dfrac{1}{3}y_1^3+\dfrac{b}{2}y_1^2+a_2y_2+\cdots+a_ny_n \nonumber
    \end{align}
for some constants $a_2,\ldots,a_n,b\in\bbR$.
\item $\Pi$, when restricted to $\Sigma$, is a self-transverse immersion.
\end{itemize}
\end{defin}

\begin{rem}
Note, if $n>1$, $L$ having simple front singularities is not a generic condition. However, if $n>2$, there is an $h$-principle which guarantees that, unless there is a homotopical reason for singularities of codimension larger than 1 to exist, we can arrange $L$ to have simple front singularities via Legendrian isotopy, cf. \cite{Ent99,AG17}. Meanwhile, if $n=2$, we can arrange $L$ to have only standard cusp-edge singularities and swallowtail singularities. The results of \cite{Ekh07} still hold in this case; therefore, the results of the present article should also remain unaffected.
\end{rem}

We will refer to $\Pi(\Sigma)\subset M$ as the \emph{singular set}. Via a generic perturbation, we may assume $\Pi\vert_\Sigma$ is self-transverse in the following sense: there is a stratification 
    \begin{equation}
    \Pi(\Sigma)=\Sigma_1\supset\cdots\supset\Sigma_K,
    \end{equation}
where $\Sigma_i$ is the set of self-intersection points of $\Pi\vert_\Sigma$ of multiplicity at least $i$, such that $\Sigma_i\subset M$ is codimension $i$ (hence, $K\leq n$). 

Locally, we may describe $L$ as the multi 1-jet graph of \emph{local functions} defined on $M$, as follows. 
\begin{enumerate}
\item First, if $p\in M-\Pi(\Sigma)$, then there exists an open neighborhood $U\subset M-\Pi(\Sigma)$ of $p$ and $k$ disjoint open sets $U_1,\ldots,U_k\subset L-\Sigma$ such that 
    \begin{equation}
    U_i\equiv\big\{(x,d(f_i)_x,f_i(x):x\in U\big\}, \;\; f_i\in C^\infty(U)
    \end{equation}
and $\Pi\vert_{U_i}:U_i\to U$ is a diffeomorphism. Each $U_i$ is referred to as a \emph{smooth sheet} lying over $p$.
\item Second, a point $p\in\Sigma_j-\Sigma_{j+1}$ may have smooth sheets lying over it. In addition, there exists an open neighborhood $V\subset\Sigma_j$ of $p$ and $j$ disjoint open sets $V_1,\ldots,V_j\subset\Sigma$ such that $\Pi\vert_{V_i}:V_i\to V$ is a diffeomorphism. If we consider a sufficiently small neighborhood $U\subset M$ of $p$, then each $\Pi(V_i)$ subdivides $U$ into two components $\{U_i^\pm\}$. Moreover, there exists two functions $\{f_i^1,f_i^2\}$ defined on $U_i^+$ such that both functions: extend to $\overline{U}_i^+$, agree on $\partial U_i^+$, and have differentials whose limits agree on $\partial U_i^+$ when approached from $U_i^+$. Since $L$ has simple front singularities, we may use the coordinates \eqref{eqn:cuspedgecoordinates} near the cusp-edge to see these local functions are of the form 
    \begin{align}
    f_i^1(x)&=\dfrac{1}{3}(2x_1)^{3/2}+bx_1+a_2x_2+\cdots+a_nx_n \\
    f_i^2(x)&=-\dfrac{1}{3}(2x_1)^{3/2}+bx_1+a_2x_2+\cdots+a_nx_n,
    \end{align}
where: these functions are defined on the subset $\{x_1>0\}$, the subset $\{x_1=0\}$ corresponds to $\Pi(\Sigma)$, and the common limit of the differentials along $\{x_1=0\}$ is 
    \begin{equation}
    bdx_1+a_2x_2+\cdots+a_nx_n.
    \end{equation}
Each $V_i$ is referred to as a \emph{singular sheet} lying over $p$.
\end{enumerate}

\begin{rem}
By compactness of $M$ and $L$, there are only finitely many local functions being considered.
\end{rem}

Let $f_i$ and $f_j$ be two local functions whose domains have non-empty intersection $U\subset M$ and consider the \emph{local function difference} $f_i-f_j\in C^\infty(U)$. We call a curve $\phi:I\to M$, where $I\subset\bbR$ is an interval, a \emph{gradient flow} of $L$ if $\phi$ is a solution to the gradient flow equation 
    \begin{equation}
    \partial_s\phi(s)=-\grad(f_i-f_j)\big(\phi(s)\big)
    \end{equation}
for some local function difference $f_i-f_j$. The \emph{1-jet lift} of $\phi$, where $\phi$ is a gradient flow of $L$ associated to $f_i-f_j$, is the ordered pair 
    \begin{equation}
    \Big(\widetilde{\phi}^1,\widetilde{\phi}^2\Big)
    \end{equation}
of continuous lifts $\widetilde{\phi}^k:I\to L$ of $\phi$ such that (1) $\widetilde{\phi}^1$ resp. $\widetilde{\phi}^2$ lies in the sheet determined by $f_i$ resp. $f_j$ and (2) the equality
    \begin{equation}
    \widetilde{\phi}^1(t_0)=\widetilde{\phi}^2(s_0)
    \end{equation}
holds if and only if $t_0=s_0$ and the point $\widetilde{\phi}^1(t_0)=\widetilde{\phi}^2(s_0)$ is a singular point lying over $\phi(t_0)=\phi(s_0)$ (i.e., the two components of the 1-jet lift may only meet in $\Sigma$). Analogously, we define the \emph{cotangent lift} of $\phi$ to be the ordered pair 
    \begin{equation}
    \Big(\overline{\phi}^1\equiv\Pi_\bbC\circ\widetilde{\phi}^1,\overline{\phi}^2\equiv\Pi_\bbC\circ\widetilde{\phi}^2\Big).
    \end{equation}

\begin{defin}\label{defin:floworientation}
The \emph{flow orientation} of $\widetilde{\phi}^1$ at $\widetilde{p}\in L$ is given by the unique lift of the vector 
    \begin{equation}
    -\grad(f_i-f_j)\big(\Pi(\widetilde{p})\big)\in T_{\Pi(\widetilde{p})}M
    \end{equation}
to $T_{\widetilde{p}}L$. The \emph{flow orientation} of $\widetilde{\phi}^2$ at $\widetilde{p}\in L$ is given by the unique lift of the vector 
    \begin{equation}
    -\grad(f_j-f_i)\big(\Pi(\widetilde{p})\big)\in T_{\Pi(\widetilde{p})}M
    \end{equation}
to $T_{\widetilde{p}}L$.
\end{defin}

\begin{rem}
Definition \ref{defin:floworientation} induces a flow orientation on the cotangent lift.
\end{rem}

By existence and uniqueness theorems of ODEs, we observe that any gradient flow of $L$ has a maximal interval of definition. Any gradient flow of $L$ defined on its maximal interval of definition will be referred to as \emph{maximally extended}.

\begin{lem}[Lemma 2.8 in \cite{Ekh07}]
Let $\Phi:I\to M$ be a maximally extended gradient flow of $L$ associated to $f_i-f_j$. 
\begin{enumerate}
\item If $I$ is an interval with an infinite end (i.e., $I=(a,+\infty)$ or $I=(-\infty,b)$ with $a,b\in\bbR\cup\{\pm\infty\}$), then 
    \begin{equation}
    \lim_{t\to\pm\infty}\Phi(t)\in\crit(f_i-f_j),
    \end{equation}
where $\crit(f_i-f_j)$ is the set of critical points of $f_i-f_j$. 
\item If $I$ is an interval with a compact end (i.e., $I=(a,\ell]$ or $I=[\ell,b)$ with $a,b\in\bbR\cup\{\pm\infty\}$ and $\ell\in\bbR$), then 
    \begin{equation}
    \lim_{t\to\ell}\widetilde{\Phi}^k(t)\in\Sigma
    \end{equation}
for $k=1$ or $k=2$. In particular, $\Phi(\ell)$ lies in $\Pi(\Sigma)$.
\end{enumerate}
\end{lem}

\begin{defin}
Let $\phi:I\to M$ be a gradient flow of $L$ associated to $f_i-f_j$ and $\Phi$ its maximal extension. We say $\phi$ is: 
\begin{itemize}
\item a \emph{Morse flow} if $\Phi$ connects two points of $\crit(f_i-f_j)$, 
\item a \emph{fold terminating flow} if $\Phi$ connects $\crit(f_i-f_j)$ to $\Pi(\Sigma)$,
\item a \emph{fold emanating flow} if $\Phi$ connects $\Pi(\Sigma)$ to $\crit(f_i-f_j)$,
\item or a \emph{singular flow} if $\Phi$ connects two points of $\Pi(\Sigma)$.
\end{itemize}
\end{defin}

\subsection{Morse flow trees}
Let 
    \begin{equation}
    T=\big(\vertex(T),\edgef(T)\big),
    \end{equation}
be a finite tree. We say that $T$ is a \emph{source tree} if, for each $v\in\vertex(T)$, we have a cyclic ordering of the set 
    \begin{equation}
    \{e\in\edgef(T):v\in e\}.
    \end{equation}
    
\begin{defin}[Definition 2.10 in \cite{Ekh07}]\label{defin:morseflowtree}
A \emph{Morse flow tree} of $L$ is a continuous map $\phi:T\to M$, where $T$ is a source tree, that satisfies the following conditions.
\begin{enumerate}
\item For each $e\in\edgef(T)$, $\phi:e\to M$ is an injective parameterization of a gradient flow of $L$. (In particular, $\phi:e\to M$ is non-constant.)
\item Let $v\in\vertex(T)$ be a $k$-valent vertex with cyclically ordered edges $e_1,\ldots,e_k$. Recall, associated to $e_j$ is its cotangent lift 
    \begin{equation}
    \Big(\overline{\phi}^1_j,\overline{\phi}^2_j\Big).
    \end{equation}
We denote by $\overline{\phi}^{in}_j$ resp. $\overline{\phi}^{out}_j$ the component of $\Big(\overline{\phi}^1_j,\overline{\phi}^2_j\Big)$ which points towards resp. away from $v$. We require the following equality: 
    \begin{equation}
    \overline{p}_j\equiv\overline{\phi}^{in}_j(v)=\overline{\phi}^{out}_{j+1}(v)\in\Pi_\bbC(L).
    \end{equation}
\item We require that all of the cotangent lifts of all of the edges of $T$ fit together to give an oriented loop in $\Pi_\bbC(L)$.
\end{enumerate}
\end{defin}

\begin{rem}
Observe, we have the following kinds of 2-valent vertices allowed in property (2) of Definition \ref{defin:morseflowtree}. It may be the case that $v$ is a 2-valent vertex such that the gradient flow $\phi_1$ of $L$ associated to $e_1$ and the gradient flow $\phi_2$ of $L$ associated to $e_2$ are both gradient flows of $L$ associated to $f_i-f_j$ and $v\notin\crit(f_i-f_j)$. We call these kinds of 2-valent vertices \emph{removable}. Removable vertices are necessary for (pre)-Floer-Gromov convergence as they can be approached by collapsing edges in a sequence of (pre)stable Morse flow trees, cf. Figure \ref{fig:ekholm19} where a removable vertex is approached by the collision of a $Y_0$-vertex with an adjacent switch and end (see Table \ref{table:vertices} for the definition of these vertices).
\end{rem}

For brevity, we will usually denote Morse flow trees by $T$, i.e., we omit the parameterization map. Also, we will denote by $\widetilde{T}$ resp. $\overline{T}$ the 1-jet resp. cotangent lift of $T$ defined by lifting each edge. 

\begin{defin}\label{defin:node1}
Let $v\in\vertex(T)$ be a vertex of a Morse flow tree $T$ and consider a non-trivial decomposition
     \begin{equation}
    I\amalg J=\{e\in\edgef(T):v\in e\},
    \end{equation}
where both $I$ and $J$ consist of consecutively ordered edges. We say $v$ has a \emph{node} if there exists a decomposition 
    \begin{equation}
    T_I\cup_{v_I\sim v_J} T_J=T,
    \end{equation}
where $I$ and $J$ are as above, such that $T_I$ resp. $T_J$ is a Morse flow tree containing a vertex $v_I$ resp. $v_J$ satisfying
    \begin{equation}
    I=\{e\in\edgef(T_I):v_I\in e\}\;\;\mathrm{resp.}\;\;J=\{e\in\edgef(T_I):v_J\in e\}.
    \end{equation}
We say $T$ is \emph{non-nodal} if it contains no nodes.
\end{defin}

Observe, in order for $T$ to have a node at a $k$-valent vertex $v$, it must be the case that there exists distinct edges $e_i\neq e_j$, with $1\leq i<j\leq k$, such that
    \begin{equation}
    \overline{\phi}^{in}_i(v)=\overline{\phi}^{out}_j(v).
    \end{equation}
In this case, we may then split off the edges $e_i,e_{i+1},\ldots,e_{j-1},e_j$ from $T$. By inductively splitting off edges, vertex by vertex in $T$, we have the following lemma (cf. \cite[Lemma 2.11]{Ekh07}).

\begin{lem}
There exists a decomposition of $T$ into non-nodal Morse flow trees:
    \begin{equation}
    T=T_{I_1}\cup_{v_{I_1}\sim v_{I_2}'}T_{I_2}\cup_{v_{I_2}\sim v_{I_2}'}\cup\cdots\cup_{v_{I_{s-1}}\sim v_{I_s}'} T_{I_s}.
    \end{equation}
Moreover, this decomposition is unique up to relabeling/reordering.
\end{lem}

\begin{rem}\label{rem:node1}
A simple example of a node is given by collapsing the edge between two adjacent $Y_0$-vertices.
\end{rem}

\begin{table}[htbp]
\begin{center}
\caption{{\bf A complete list of vertices occurring in rigid Morse flow trees (borrowed from \cite[Table 1]{Kar24})}}\label{table:vertices}
\tabulinesep=1.15mm
\begin{tabu}{>{\centering\arraybackslash}m{10mm} | >{\arraybackslash}m{60mm}|>{\centering\arraybackslash}m{45mm}}
    \hline Type  & Description & Local picture of front \\
    \hline $P_{1,+}$ & \emph{Positive 1-valent punctures}, not contained in $\Pi(\Sigma)$. & \includegraphics[scale=.5]{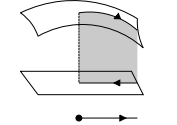}\\
    \hline $P_{1,-}$ &  \emph{Negative 1-valent punctures}, not contained in $\Pi(\Sigma)$. & \includegraphics[scale=.5]{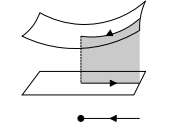}\\
   \hline $P_{2,+}$ & \emph{Positive 2-valent punctures $p$}, not contained in  $\Pi(\Sigma)$, with index $I(p)=0$. & \includegraphics[scale=.5]{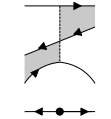}\\
   \hline $P_{2,-}$ & \emph{Negative 2-valent punctures $p$}, not contained in  $\Pi(\Sigma)$, with index $I(p)=n$. & \includegraphics[scale=.5]{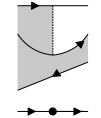}\\
   \hline $Y_0$ &  \emph{3-valent $Y_0$-vertices}, not contained in  $\Pi(\Sigma)$ and not containing any punctures. & \includegraphics[scale=.5]{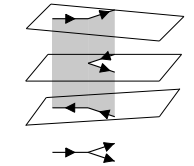}\\
   \hline end & \emph{1-valent end-vertices $e$},  contained in $\Pi(\Sigma)$, not containing any punctures, meeting $\Pi(\Sigma)$ transversely, the $1$-jet lift of $T$ through $e$ is traversing $\Sigma$ downwards. & \includegraphics[scale=.5]{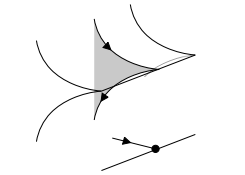}\\
   \hline switch & \emph{2-valent switch-vertices $s$}, contained in  $\Pi(\Sigma)$,  not containing any punctures, tangent to $\Pi(\Sigma)$,  one of the $1$-jet lifts of $T$ through $v$ is traversing $\Sigma$ upwards while the other 1-jet lift of $T$ through $s$ is not contained in $\Sigma$. & \includegraphics[scale=.5]{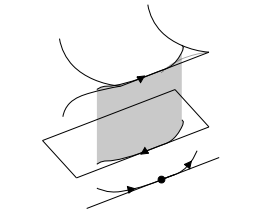}\\
   \hline $Y_1$ & \emph{3-valent $Y_1$-vertices  $v$},  contained in $\Pi(\Sigma)$, not containing any punctures, meeting $\Pi(\Sigma)$ transversely, one of the $1$-jet lifts of $T$ through $v$ is traversing $\Sigma$ upwards while the other two 1-jet lifts of $T$ through $v$ are not contained in $\Sigma$. & \includegraphics[scale=.35]{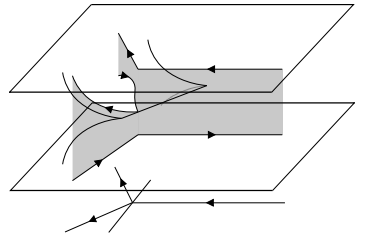}\\
   \hline
\end{tabu}
\end{center}
\end{table}

Again, let $v$ be a $k$-valent vertex, with cyclically ordered edges $e_1,\ldots,e_k$, of a Morse flow tree $T$ and consider two paired cotangent lifts $\overline{\phi}^{in}_j$ and $\overline{\phi}^{out}_{j+1}$.

\begin{defin}\label{defin:puncture}
We say that $v$ contains a \emph{puncture} after $e_j$ if 
    \begin{equation}\label{eqn:puncture}
    \widetilde{\phi}^{in}_j(v)\neq\widetilde{\phi}^{out}_{j+1}(v)\in L.
    \end{equation}
\end{defin}

Observe, if \eqref{eqn:puncture} holds, then the aforementioned points must be endpoints of a Reeb chord on $L$. We say this puncture is \emph{positive} if 
    \begin{equation}
    z\Big(\widetilde{\phi}^{in}_j(v)\Big)>z\Big(\widetilde{\phi}^{out}_{j+1}(v)\Big)
    \end{equation}
and \emph{negative} otherwise. \cite[Lemma 2.11]{Ekh07} shows we may decompose any Morse flow tree into a union of Morse flow trees where each vertex contains at most 1 puncture (in particular, the supplied proof implies a vertex with more than 1 puncture is nodal). Also, Lemma 2.13 of \emph{loc. cit.} shows every Morse flow tree has at least 1 positive puncture.

Now, we will assume $L$ satisfies the \emph{preliminary transversality condition} detailed in \cite[Pages 1101-1103]{Ekh07}. We will not reproduce the condition here since its explanation is rather long-winded. In short, the preliminary transversality condition details how the various sheets of $L$ are allowed to meet each other, and the main consequence is that gradient flows of $L$ meet $\Pi(\Sigma)$ in a locally stable fashion. In particular, the main use of this condition in the present article is that it implies any gradient flow of $L$ meets $\Pi(\Sigma)$ with order of contact at most $n$ (in particular, no gradient flow of $L$ flows along $\Pi(\Sigma)$). Finally, the preliminary transversality condition may be achieved after a generic perturbation of $L$.

Let $T$ be a Morse flow tree. We denote by: 
\begin{itemize}
\item $P(T)$ the set of positive punctures of $T$, 
\item $Q(T)$ the set of negative punctures of $T$, 
\item and $R(T)$ set of vertices of $T$ which do not contain punctures. 
\end{itemize}
Let $r\in R(T)$ and consider a cusp point $x\in\Sigma$ over $r$ lying in $\widetilde{T}$; we define $\widetilde{\mu}(x)=1$ resp. $\widetilde{\mu}(x)=-1$ if the incoming arc of $\widetilde{T}$ lies in the upper resp. lower sheet of $L$ and the outgoing arc of $\widetilde{T}$ lies in the lower resp. upper sheet of $L$. In all other cases, we define $\widetilde{\mu}(x)=0$. The \emph{Maslov content} of $r$ is 
    \begin{equation}
    \mu(r)\equiv\sum\widetilde{\mu}(x),
    \end{equation}
where the sum is over all cusp points $x\in\Sigma$ over $r$ lying in $\widetilde{T}$.

\begin{defin}\label{defin:fd}
The \emph{formal dimension} of a Morse flow tree $T$ is
\begin{multline}
\dim(T)\equiv\sum_{p\in P(T)}I(p)+\sum_{q\in Q(T)}\big(n-I(q)\big)+\sum_{r\in R(T)}\mu(r)- \\
n\big(\abs{P(T)}+\abs{Q(T)}-1\big)+\big(\abs{P(T)}+\abs{Q(T)}-3\big),
\end{multline}
where $I(p)$ resp. $I(q)$ is the Morse index of $p$ resp. $q$.
\end{defin}

A Morse flow tree of formal dimension 0 is called \emph{rigid}. The vertices of rigid Morse flow trees are completely classified, cf. Table \ref{table:vertices}.

\begin{lem}[Lemma 3.12 in \cite{Ekh07}]
Let $P,D\in\bbZ_{>0}$. There exists a constant $N_{P,D}\in\bbZ_{>0}$ such that, for any Morse flow tree $T$ with at most $P$ positive punctures and of formal dimension at most $D$, the domain of $T$ has at most $N_{P,D}$ edges.
\end{lem}

\section{Compactifying edges}
Let $f_i-f_j$ be a local function difference defined on some open subset $U\subset M$. Since the Reeb chords of $L$ are in 1-1 correspondence with the self-transverse double points of $\Pi_\bbC(L)$, and $L$ is chord generic, it follows that $f_i-f_j\in C^\infty(U)$ is a Morse function whose set of critical points $\crit(f_i-f_j)$ is finite. Moreover, by a generic perturbation of $L$, we may assume no critical point of any local function difference lies in $\Pi(\Sigma)$. Floer-Gromov convergence will essentially reduce to a statement about convergence of each associated sequence of edges, therefore, in this section, we will compactify the various moduli spaces of gradient flows of our various local function differences via standard techniques of Morse theory essentially following \cite{AD14}.

In order to simplify various technical arguments, we will make the following assumption.

\begin{assu}\label{assu:main}
The Riemannian metric $g$ on $M$ is Euclidean in any Morse neighborhood of any critical point of any local function difference. In particular, in a Morse neighborhood of any critical point of $f_i-f_j$, $\grad(f_i-f_j)$ restricts to the usual gradient with respect to the usual Euclidean metric.
\end{assu}

\begin{rem}
Recall, a Morse neighborhood $\Omega(c)\subset U$ of $c\in\crit(f_i-f_j)$ is an open subset such that $f_i-f_j$ takes the form
    \begin{equation}
    (f_i-f_j)(c)-x_1^2-\cdots-x_{I(c)}^2+x_{I(c)+1}^2+\cdots+x_n^2,
    \end{equation}
where $I(c)$ is the Morse index of $c$.
\end{rem}

\begin{rem}
Of course, Assumption \ref{assu:main} cannot be generically achieved. However, given $f_i-f_j$, we can always find a \emph{gradient-like} vector field $X_{ij}$ for $f_i-f_j$, i.e., a vector field $X_{ij}$ on $U$ satisfying:
\begin{enumerate}
\item $X_{ij}(f_i-f_j)(x)\leq0$ with equality if and only if $x\in\crit(f_i-f_j)$;
\item and, on any $\Omega(c)$ with $c\in\crit(f_i-f_j)$, $X_{ij}$ restricts to the usual gradient of $f_i-f_j$ with respect to the usual Euclidean metric.
\end{enumerate}
Moreover, we can find a gradient-like vector field whose flow is sufficiently close to the gradient flow of $f_i-f_j$, so there is really no harm in using Assumption \ref{assu:main}.
\end{rem}

\subsection{Morse theory preliminaries}
Using Assumption \ref{assu:main}, we will describe a convenient form that our Morse neighborhoods now take, cf. \cite[Pages 24-29]{AD14}. Let $c\in\crit(f_i-f_j)$ and consider a Morse neighborhood
    \begin{equation}
    \varphi:\Omega(c)\xrightarrow{\sim}V\subset\bbR^n
    \end{equation}
of $c$, where we assume $\varphi(c)=0$. We will denote by $Q$ the Hessian of $(f_i-f_j)\circ\varphi^{-1}$ at the origin: 
    \begin{equation}
    Q\equiv d^2\big((f_i-f_j)\circ\varphi^{-1}\big)_0.
    \end{equation}
We have that $Q$ is negative-definite resp. positive-definite on a $I(c)$-dimensional subspace $V_-\subset\bbR^n$ resp. $\big(n-I(c)\big)$-dimensional subspace $V_+\subset\bbR^n$ such that $V_-$ and $V_+$ are complimentary. In particular, it follows that $V$ is of the form 
    \begin{equation}
    V=\big\{x\in\bbR^n:-\epsilon<Q(x)<\epsilon,\abs{x_-}^2\abs{x_+}^2\leq\eta(\epsilon+\eta)\big\},
    \end{equation}
where $\epsilon,\eta\in\bbR_{>0}$ are sufficiently small and 
    \begin{equation}
    x\equiv x_-\oplus x_+\in V_-\oplus V_+.
    \end{equation}
Moreover, the boundary of $V$ consists of three parts: 
    \begin{align}
    \partial_\pm V&\equiv\big\{x\in\bbR^n:Q(x)=\pm\epsilon,\abs{x_\mp}^2\leq\eta\big\}, \\
    \partial_0V&\equiv\big\{x\in\bbR^n:\abs{x_-}^2\abs{x_+}^2=\eta(\epsilon+\eta)\big\}.
    \end{align}
Recall, the \emph{unstable} and \emph{stable} manifolds of $c$ are defined to be
    \begin{align}
    W^u(c;f_i-f_j)&=\big\{p\in U:\lim_{s\to-\infty}\phi_p(s)=c\big\}, \\
    W^s(c;f_i-f_j)&=\big\{p\in U:\lim_{s\to+\infty}\phi_p(s)=c\big\},
    \end{align}
where $\phi_p$ is the gradient flow of $f_i-f_j$ determined by $p$, respectively. (We will omit $f_i-f_j$ from the notation when it is clear from context). We have that $W^u(c)$ resp. $W^s(c)$ is diffeomorphic to a disk $D^{I(c)}$ resp. $D^{n-I(c)}$. Moreover, by the description $\Omega(c)$ just described, we see that 
    \begin{align}
    \partial_-V&\cong W^u(c)\cap\partial\Omega(c)\cong S^{I(c)-1}, \\
    \partial_+V&\cong W^s(c)\cap\partial\Omega(c)\cong S^{n-I(c)-1}.
    \end{align}

We will require the following two lemmata in the sequel.

\begin{lem}[Lemma 3.2.5 in \cite{AD14}]\label{lem:prelim1}
Let $x\in U-\crit(f_i-f_j)$ and consider a sequence $\{x_\nu\}\subset U$ converging to $x$. Suppose $\{y_\nu\}\subset U$ is a sequence of points resp. $y\in U$ is a point lying on the same gradient flow(s) of $f_i-f_j$ as $\{x_\nu\}$ resp. $x$, and 
    \begin{equation}
    (f_i-f_j)(y_\nu)=(f_i-f_j)(y),\;\;\forall\nu,
    \end{equation}
then $\{y_\nu\}$ converges to $y$.
\end{lem}

\begin{lem}\label{lem:prelim2}
Consider a sequence of maximally extended gradient flows $\{\Phi_\nu\}$ resp. a maximally extended fold terminating gradient flow $\Phi$ of $f_i-f_j$. Suppose there exists a sequence $\{y_\nu\}$ of points lying on $\{\Phi_\nu\}$ converging to the ending point $y$ of $\Phi$, then $\{\Phi_\nu\}$ is actually a sequence of maximally extended fold terminating gradient flows whose associated sequence of ending points converges to $y$.
\end{lem}

\begin{proof}
This is immediate after investigating $\grad(f_i-f_j)$ in a chart near $y$.
\end{proof}

\subsection{Edge convergence}
\subsubsection{Notation}
In this part, we will set notation. Let $a,b\in\crit(f_i-f_j)$, where $a$ and $b$ are not necessarily distinct. We may define the following moduli spaces: 
\begin{itemize}
\item the moduli space of parameterized Morse flows of $f_i-f_j$ connecting $a$ to $b$, denoted
    \begin{equation}
    \scrM_0(a,b;f_i-f_j);
    \end{equation}
\item the moduli space of parameterized fold terminating flows of $f_i-f_j$ starting at $a$, denoted
    \begin{equation}
    \scrM_0\big(a,\Pi(\Sigma);f_i-f_j\big);
    \end{equation}
\item the moduli space of parameterized fold emanating flows of $f_i-f_j$ ending at $b$, denoted
    \begin{equation}
    \scrM_0\big(\Pi(\Sigma),b;f_i-f_j\big);
    \end{equation}
\item the moduli space of parameterized singular flows of $f_i-f_j$, denoted
    \begin{equation}
    \scrM_0\big(\Pi(\Sigma);f_i-f_j\big).
    \end{equation}
\end{itemize}
We will omit the subscript in the notation to define the analogous moduli space consisting only of maximally extend parameterized gradient flows of $f_i-f_j$ of the appropriate type. When $a\neq b$, there is a smooth proper $\bbR$-action on each moduli space given by time-shift: 
    \begin{equation}
    \big(r_0,\phi(s)\big)\mapsto\phi(s+r_0),\;\;r\in\bbR.
    \end{equation}
In particular, we may define the following moduli spaces: 
\begin{itemize}
\item the moduli space of Morse flows of $f_i-f_j$ connecting $a$ to $b$, denoted
    \begin{equation}
    \scrL_0(a,b;f_i-f_j)\equiv\scrM_0(a,b;f_i-f_j)/\bbR;
    \end{equation}
\item the moduli space of fold terminating flows of $f_i-f_j$ starting at $a$, denoted
    \begin{equation}
    \scrL_0\big(a,\Pi(\Sigma);f_i-f_j\big)\equiv\scrM_0\big(a,\Pi(\Sigma);f_i-f_j\big)/\bbR;
    \end{equation}
\item the moduli space of fold emanating flows of $f_i-f_j$ ending at $b$, denoted
    \begin{equation}
    \scrL_0\big(\Pi(\Sigma),b;f_i-f_j\big)\equiv\scrM_0\big(\Pi(\Sigma),b;f_i-f_j\big)/\bbR;
    \end{equation}
\item the moduli space of singular flows of $f_i-f_j$, denoted
    \begin{equation}
    \scrL_0\big(\Pi(\Sigma);f_i-f_j\big)\equiv\scrM_0\big(\Pi(\Sigma);f_i-f_j\big)/\bbR.
    \end{equation}
\end{itemize}
Again, we will omit the subscript in the notation to define the analogous moduli space consisting only of maximally extend gradient flows of $f_i-f_j$ of the appropriate type. Let $c_1,\ldots,c_q\in\crit(f_i-f_j)$. We also let 
    \begin{equation}
    c_0=a,\;\;c_{q+1}=b.
    \end{equation}

\begin{defin}
For $q\in\bbZ_{>0}$, a \emph{$q$-times broken Morse flow} of $f_i-f_j$ connecting $a$ to $b$, denoted
    \begin{equation}
    \lambda\equiv(\lambda^1,\ldots,\lambda^{q+1}),
    \end{equation}
is a concatenation of gradient flows of $f_i-f_j$: 
    \begin{equation}
    \lambda^1\in\scrL_0(c_0,c_1;f_i-f_j),\;\;\lambda^k\in\scrL(c_{k-1},c_k;f_i-f_j),\;\;\lambda^{q+1}\in\scrL_0(c_q,c_{q+1};f_i-f_j),
    \end{equation}
where $2\leq k\leq q$. A \emph{0-times broken Morse flow} of $f_i-f_j$ connecting $a$ to $b$ is simply a Morse flow. We denote by 
    \begin{equation}
    \scrL_{q,0}(a,b;f_i-f_j)
    \end{equation}
the moduli space of broken Morse flows of $f_i-f_j$ connecting $a$ to $b$.
\end{defin}

\begin{defin}
For $q\in\bbZ_{>0}$, a \emph{$q$-times broken fold terminating flow} of $f_i-f_j$ starting at $a$, denoted
    \begin{equation}
    \lambda\equiv(\lambda^1,\ldots,\lambda^{q+1}),
    \end{equation}
is a concatenation of gradient flows of $f_i-f_j$: 
    \begin{equation}
    \lambda^1\in\scrL_0(c_0,c_1;f_i-f_j),\;\;\lambda^k\in\scrL(c_{k-1},c_k;f_i-f_j),\;\;\lambda^{q+1}\in\scrL_0\big(c_q,\Pi(\Sigma);f_i-f_j\big),
    \end{equation}
where $2\leq k\leq q$. A \emph{0-times broken fold terminating} of $f_i-f_j$ starting at $a$ is simply a fold terminating flow. We denote by 
    \begin{equation}
    \scrL_{q,0}\big(a,\Pi(\Sigma);f_i-f_j\big)
    \end{equation}
the moduli space of fold terminating flows of $f_i-f_j$ starting at $a$.
\end{defin}

\begin{defin}
For $q\in\bbZ_{>0}$, a \emph{$q$-times broken fold emanating flow} of $f_i-f_j$ ending at $b$, denoted
    \begin{equation}
    \lambda\equiv(\lambda^1,\ldots,\lambda^{q+1}),
    \end{equation}
is a concatenation of gradient flows of $f_i-f_j$: 
    \begin{equation}
    \lambda^1\in\scrL_0\big(\Pi(\Sigma),c_1;f_i-f_j\big),\;\;\lambda^k\in\scrL(c_{k-1},c_k;f_i-f_j),\;\;\lambda^{q+1}\in\scrL_0(c_q,c_{q+1};f_i-f_j),
    \end{equation}
where $2\leq k\leq q$. A \emph{0-times broken fold emanating} of $f_i-f_j$ ending at $b$ is simply a fold terminating flow. We denote by 
    \begin{equation}
    \scrL_{q,0}\big(\Pi(\Sigma),b;f_i-f_j\big)
    \end{equation}
the moduli space of fold emanating flows of $f_i-f_j$ ending at $b$.
\end{defin}

\begin{defin}
For $q\in\bbZ_{>0}$, a \emph{$q$-times broken singular flow} of $f_i-f_j$, denoted
    \begin{equation}
    \lambda\equiv(\lambda^1,\ldots,\lambda^{q+1}),
    \end{equation}
is a concatenation of gradient flows of $f_i-f_j$: 
    \begin{equation}
    \lambda^1\in\scrL_0\big(\Pi(\Sigma),c_1;f_i-f_j\big),\;\;\lambda^k\in\scrL(c_{k-1},c_k;f_i-f_j),\;\;\lambda^{q+1}\in\scrL_0\big(c_q,\Pi(\Sigma);f_i-f_j\big),
    \end{equation}
where $2\leq k\leq q$. A \emph{0-times broken singular flow} of $f_i-f_j$ is simply a singular flow. We denote by 
    \begin{equation}
    \scrL_{q,0}\big(\Pi(\Sigma);f_i-f_j\big)
    \end{equation}
the moduli space of singular flows of $f_i-f_j$.
\end{defin}

Analogously, we may define a maximally extended $q$-times broken gradient flow of $f_i-f_j$ of the appropriate type by restricting to the case in the previous four definitions that the first and last component of the considered broken gradient flow are maximally extended; we omit the subscript in the notation to define the moduli space of maximally extended $q$-times broken gradient flows of $f_i-f_j$ of the appropriate type. We may now define the following stratified moduli spaces: 
\begin{align}
\overline{\scrL}_0(a,b;f_i-f_j)&\equiv\bigcup_{q\geq0}\scrL_{q,0}(a,b;f_i-f_j), \\
\overline{\scrL}_0\big(a,\Pi(\Sigma);f_i-f_j\big)&\equiv\bigcup_{q\geq0}\scrL_{q,0}\big(a,\Pi(\Sigma);f_i-f_j\big), \\
\overline{\scrL}_0\big(\Pi(\Sigma),b;f_i-f_j\big)&\equiv\bigcup_{q\geq0}\scrL_{q,0}\big(\Pi(\Sigma),b;f_i-f_j\big), \\
\overline{\scrL}_0\big(\Pi(\Sigma);f_i-f_j\big)&\equiv\bigcup_{q\geq0}\scrL_{q,0}\big(\Pi(\Sigma);f_i-f_j\big).
\end{align}

\begin{rem}
Of course, there are only finitely many non-empty sets in the unions of the previous four moduli spaces since (1) $\crit(f_i-f_j)$ is a finite set and (2) $f_i-f_j$ strictly decreases along any gradient flow.
\end{rem}

Again, we will omit the subscript in the notation to define the analogous stratified moduli spaces consisting of only maximally extended objects. Finally, we define
    \begin{multline}
    \overline{\scrL}_0(f_i-f_j)\equiv\overline{\scrL}_0\big(\Pi(\Sigma);f_i-f_j)\cup \\
    \bigcup_{a,b\in\crit(f_i-f_j)}\Big\{\overline{\scrL}_0(a,b;f_i-f_j)\cup\overline{\scrL}_0\big(a,\Pi(\Sigma);f_i-f_j\big)\cup\overline{\scrL}_0\big(\Pi(\Sigma),b;f_i-f_j\big)\Big\},
    \end{multline}
and its analogous counterpart given by dropping the subscript in the notation.

\subsubsection{Topologies}
In this part, we will define a topology on each moduli space:
    \begin{equation}
    \overline{\scrL}_0(a,b;f_i-f_j),\;\;\overline{\scrL}_0\big(a,\Pi(\Sigma);f_i-f_j\big),\;\;\overline{\scrL}_0\big(\Pi(\Sigma),b;f_i-f_j\big),\;\;\overline{\scrL}_0\big(\Pi(\Sigma);f_i-f_j\big),
    \end{equation}
and their maximally extended counterparts, via constructing a subbasis; then, we will show each topology is compact.

First, we consider both $\overline{\scrL}_0(a,b;f_i-f_j)$ and $\overline{\scrL}(a,b;f_i-f_j)$. Again, let 
    \begin{equation}
    \lambda=(\lambda^1,\ldots,\lambda^{q+1})\in\overline{\scrL}_0(a,b;f_i-f_j)
    \end{equation}
connect $a=c_0$ to $c_1$, $c_1$ to $c_2$, ..., and $c_q$ to $b=c_{q+1}$. We denote by 
    \begin{equation}
    \Lambda\equiv(\Lambda^1,\ldots,\Lambda^{q+1})\in\overline{\scrL}(a,b;f_i-f_j)
    \end{equation}
the maximal extension of $\lambda$.\footnote{Of course, $\lambda^k=\Lambda^k$ for $2\leq k\leq q$ by definition.} We fix a Morse neighborhood $\Omega(c_k)$ of each $c_k$. Let $\scrU^-_{k-1}\subset\partial\Omega(c_{k-1})$ resp. $\scrU^+_k\subset\partial\Omega(c_k)$ be a neighborhood of the exit resp. entry point of $\Lambda^k$ from $\Omega(c_{k-1})$ resp. into $\Omega(c_k)$ such that this neighborhood is contained in a level set of $f_i-f_j$. Note, such choices of $\scrU^\pm_k$'s are guaranteed by our earlier description of $\Omega(c_k)$. Also, let $\scrV^-\subset U$ resp. $\scrV^+\subset U$ be a neighborhood of the starting point of $\lambda^1$ resp. the ending point of $\lambda^{q+1}$. We define the collections
    \begin{equation}
    \scrU^-\equiv\{\scrU^-_0,\ldots,\scrU^-_q,\scrV^-\}\;\;\textrm{and}\;\;\scrU^+\equiv\{\scrU^+_1,\ldots,\scrU^+_{q+1},\scrV^+\}.
    \end{equation}
Let $r\in\bbZ_{\geq0}$ with $r+1\leq q+1$. 
\begin{enumerate}
\item We say that 
    \begin{equation}
    \Theta=(\Theta^1,\ldots,\Theta^{r+1})\in\scrW(\Lambda,\scrU^-,\scrU^+)\subset\overline{\scrL}(a,b;f_i-f_j)
    \end{equation}
if there exists non-negative integers
    \begin{equation}
    k_0=0<k_1<\cdots<k_r<k_{r+1}=q+1
    \end{equation}
such that the following conditions hold.
\begin{itemize}
\item $\Theta^s\in\scrL(c_{k_{s-1}},c_{k_s};f_i-f_j)$.
\item $\Theta^s$ exits $\Omega(c_{k_{s-1}})$ through $\scrU^-_{k_{s-1}}$ and enters $\Omega(c_{k_s})$ through $\scrU^+_{k_s}$.
\end{itemize}
\item We say that 
    \begin{equation}
    \theta\equiv(\theta^1,\ldots,\theta^{r+1})\in\scrW(\lambda,\scrU^-,\scrU^+)\subset\overline{\scrL}_0(a,b;f_i-f_j)
    \end{equation}
if the following conditions hold, where $\Theta$ is the maximal extension of $\theta$.
\begin{itemize}
\item $\Theta\in\scrW(\Lambda,\scrU^-,\scrU^+)$.
\item $\theta^1$'s starting point is in $\scrV^-$ and $\theta^{r+1}$'s ending point is in $\scrV^+$.
\end{itemize}
\end{enumerate}
The collection 
    \begin{equation}
    \big\{\scrW(\lambda,\scrU^-,\scrU^+)\big\}\;\;\textrm{resp.}\;\;\big\{\scrW(\Lambda,\scrU^-,\scrU^+)\big\}
    \end{equation}
provides a fundamental system of open neighborhoods for a metrizable topology on 
    \begin{equation}
    \overline{\scrL}_0(a,b;f_i-f_j)\;\;\textrm{resp.}\;\;\overline{\scrL}(a,b;f_i-f_j).
    \end{equation}

Second, we consider both $\overline{\scrL}_0\big(a,\Pi(\Sigma);f_i-f_j\big)$ and $\overline{\scrL}\big(a,\Pi(\Sigma);f_i-f_j\big)$. Again, let 
    \begin{equation}
    \lambda=(\lambda^1,\ldots,\lambda^{q+1})\in\overline{\scrL}_0\big(a,\Pi(\Sigma);f_i-f_j\big)
    \end{equation}
connect $a=c_0$ to $c_1$, $c_1$ to $c_2$, ..., and $c_q$ to a point $c_{q+1}$ in $\Pi(\Sigma)$. We denote by 
    \begin{equation}
    \Lambda\equiv(\Lambda^1,\ldots,\Lambda^{q+1})\in\overline{\scrL}_0\big(a,\Pi(\Sigma);f_i-f_j\big)
    \end{equation}
the maximal extension of $\lambda$. We fix a Morse neighborhood $\Omega(c_k)$ of each $c_k$, where $0\leq k\leq q$. Let $\scrU^-_{k-1}\subset\partial\Omega(c_{k-1})$ resp. $\scrU^+_k\subset\partial\Omega(c_k)$ be a neighborhood of the exit resp. entry point of $\Lambda^k$ from $\Omega(c_{k-1})$ resp. into $\Omega(c_k)$ such that this neighborhood is contained in a level set of $f_i-f_j$. Let $\scrU^+_{q+1}\subset\partial U$ be a neighborhood of $c_{q+1}$. Also, let $\scrV^-\subset U$ resp. $\scrV^+\subset U$ be a neighborhood of the starting point of $\lambda^1$ resp. the ending point of $\lambda^{q+1}$. We define the collections
    \begin{equation}
    \scrU^-\equiv\{\scrU^-_0,\ldots,\scrU^-_q,\scrV^-\}\;\;\textrm{and}\;\;\scrU^+\equiv\{\scrU^+_1,\ldots,\scrU^+_{q+1},\scrV^+\}.
    \end{equation}
Let $r\in\bbZ_{\geq0}$ with $r+1\leq q+1$. 
\begin{enumerate}
\item We say that 
    \begin{equation}
    \Theta=(\Theta^1,\ldots,\Theta^{r+1})\in\scrW(\Lambda,\scrU^-,\scrU^+)\subset\overline{\scrL}\big(a,\Pi(\Sigma);f_i-f_j\big)
    \end{equation}
if there exists non-negative integers
    \begin{equation}
    k_0=0<k_1<\cdots<k_r<k_{r+1}=q+1
    \end{equation}
such that the following conditions hold.
\begin{itemize}
\item $\Theta^s\in\scrL(c_{k_{s-1}},c_{k_s};f_i-f_j)$, where $1\leq s\leq r$.
\item $\Theta^s$ exits $\Omega(c_{k_{s-1}})$ through $\scrU^-_{k_{s-1}}$ and enters $\Omega(c_{k_s})$ through $\scrU^+_{k_s}$, where $1\leq s\leq r$.
\item $\Theta^{r+1}\in\scrL\big(c_{k_r},\Pi(\Sigma);f_i-f_j\big)$.
\item $\Theta^{r+1}$ exits $\Omega(c_{k_r})$ through $\scrU^-_{k_r}$ and ends at a point in $\scrU^+_{q+1}$.
\end{itemize}
\item We say that 
    \begin{equation}
    \theta\equiv(\theta^1,\ldots,\theta^{r+1})\in\scrW(\lambda,\scrU^-,\scrU^+)\subset\overline{\scrL}_0\big(a,\Pi(\Sigma);f_i-f_j\big)
    \end{equation}
if the following conditions hold, where $\Theta$ is the maximal extension of $\theta$.
\begin{itemize}
\item $\Theta\in\scrW(\Lambda,\scrU^-,\scrU^+)$.
\item $\theta^1$'s starting point is in $\scrV^-$ and $\theta^{r+1}$'s ending point is in $\scrV^+$.
\end{itemize}
\end{enumerate}
The collection 
    \begin{equation}
    \big\{\scrW(\lambda,\scrU^-,\scrU^+)\big\}\;\;\textrm{resp.}\;\;\big\{\scrW(\Lambda,\scrU^-,\scrU^+)\big\}
    \end{equation}
provides a fundamental system of open neighborhoods for a metrizable topology on 
    \begin{equation}
    \overline{\scrL}_0\big(a,\Pi(\Sigma);f_i-f_j\big)\;\;\textrm{resp.}\;\;\overline{\scrL}\big(a,\Pi(\Sigma);f_i-f_j\big).
    \end{equation}

Third, we consider both $\overline{\scrL}_0\big(\Pi(\Sigma),b;f_i-f_j\big)$ and $\overline{\scrL}\big(\Pi(\Sigma),b;f_i-f_j\big)$. Again, let 
    \begin{equation}
    \lambda=(\lambda^1,\ldots,\lambda^{q+1})\in\overline{\scrL}_0\big(\Pi(\Sigma),b;f_i-f_j\big)
    \end{equation}
connect a point $c_0$ in $\Pi(\Sigma)$ to $c_1$, $c_1$ to $c_2$, ..., and $c_q$ to $b=c_{q+1}$. We denote by 
    \begin{equation}
    \Lambda\equiv(\Lambda^1,\ldots,\Lambda^{q+1})\in\overline{\scrL}\big(\Pi(\Sigma),b;f_i-f_j\big)
    \end{equation}
the maximal extension of $\lambda$. We fix a Morse neighborhood $\Omega(c_k)$ of each $c_k$, where $1\leq k\leq q+1$. Let $\scrU^-_{k-1}\subset\partial\Omega(c_{k-1})$ resp. $\scrU^+_k\subset\partial\Omega(c_k)$ be a neighborhood of the exit resp. entry point of $\Lambda^k$ from $\Omega(c_{k-1})$ resp. into $\Omega(c_k)$ such that this neighborhood is contained in a level set of $f_i-f_j$. Let $\scrU^-_0\subset\partial U$ be a neighborhood of $c_0$. Also, let $\scrV^-\subset U$ resp. $\scrV^+\subset U$ be a neighborhood of the starting point of $\lambda^1$ resp. the ending point of $\lambda^{q+1}$. We define the collections
    \begin{equation}
    \scrU^-\equiv\{\scrU^-_0,\ldots,\scrU^-_q,\scrV^-\}\;\;\textrm{and}\;\;\scrU^+\equiv\{\scrU^+_1,\ldots,\scrU^+_{q+1},\scrV^+\}.
    \end{equation}
Let $r\in\bbZ_{\geq0}$ with $r+1\leq q+1$. 
\begin{enumerate}
\item We say that 
    \begin{equation}
    \Theta=(\Theta^1,\ldots,\Theta^{r+1})\in\scrW(\Lambda,\scrU^-,\scrU^+)\subset\overline{\scrL}\big(\Pi(\Sigma),b;f_i-f_j\big)
    \end{equation}
if there exists non-negative integers
    \begin{equation}
    k_0=0<k_1<\cdots<k_r<k_{r+1}=q+1
    \end{equation}
such that the following conditions hold.
\begin{itemize}
\item $\Theta^s\in\scrL(c_{k_{s-1}},c_{k_s};f_i-f_j)$, where $2\leq s\leq r+1$.
\item $\Theta^s$ exits $\Omega(c_{k_{s-1}})$ through $\scrU^-_{k_{s-1}}$ and enters $\Omega(c_{k_s})$ through $\scrU^+_{k_s}$, where $2\leq s\leq r+1$.
\item $\Theta^1\in\scrL\big(\Pi(\Sigma),c_1;f_i-f_j\big)$.
\item $\Theta^1$ starts at a point in $\scrU^-_0$ and enters $\Omega(c_{k_1})$ through $\scrU^+_{k_1}$.
\end{itemize}
\item We say that 
    \begin{equation}
    \theta\equiv(\theta^1,\ldots,\theta^{r+1})\in\scrW(\lambda,\scrU^-,\scrU^+)\subset\overline{\scrL}_0\big(\Pi(\Sigma),b;f_i-f_j\big)
    \end{equation}
if the following conditions hold, where $\Theta$ is the maximal extension of $\theta$.
\begin{itemize}
\item $\Theta\in\scrW(\Lambda,\scrU^-,\scrU^+)$.
\item $\theta^1$'s starting point is in $\scrV^-$ and $\theta^{r+1}$'s ending point is in $\scrV^+$.
\end{itemize}
\end{enumerate}
The collection 
    \begin{equation}
    \big\{\scrW(\lambda,\scrU^-,\scrU^+)\big\}\;\;\textrm{resp.}\;\;\big\{\scrW(\Lambda,\scrU^-,\scrU^+)\big\}
    \end{equation}
provides a fundamental system of open neighborhoods for a metrizable topology on    
    \begin{equation}
    \overline{\scrL}_0\big(\Pi(\Sigma),b;f_i-f_j\big)\;\;\textrm{resp.}\;\;\overline{\scrL}\big(\Pi(\Sigma),b;f_i-f_j\big).
    \end{equation}

Fourth, we consider both $\overline{\scrL}_0\big(\Pi(\Sigma);f_i-f_j\big)$ and $\overline{\scrL}\big(\Pi(\Sigma);f_i-f_j\big)$. Again, let 
    \begin{equation}
    \lambda=(\lambda^1,\ldots,\lambda^{q+1})\in\overline{\scrL}_0\big(\Pi(\Sigma);f_i-f_j\big)
    \end{equation}
connect a point $c_0$ in $\Pi(\Sigma)$ to $c_1$, $c_1$ to $c_2$, ..., and $c_q$ to a point $c_{q+1}$ in $\Pi(\Sigma)$. We denote by 
    \begin{equation}
    \Lambda\equiv(\Lambda^1,\ldots,\Lambda^{q+1})\in\overline{\scrL}\big(\Pi(\Sigma);f_i-f_j\big)
    \end{equation}
the maximal extension of $\lambda$. We fix a Morse neighborhood $\Omega(c_k)$ of each $c_k$, where $1\leq k\leq q$. Let $\scrU^-_{k-1}\subset\partial\Omega(c_{k-1})$ resp. $\scrU^+_k\subset\partial\Omega(c_k)$ be a neighborhood of the exit resp. entry point of $\Lambda^k$ from $\Omega(c_{k-1})$ resp. into $\Omega(c_k)$ such that this neighborhood is contained in a level set of $f_i-f_j$. Let $\scrU^-_0\subset\partial U$ resp. $\scrU^+_{q+1}$ be a neighborhood of $c_0$ resp. $c_{q+1}$. Also, let $\scrV^-\subset U$ resp. $\scrV^+\subset U$ be a neighborhood of the starting point of $\lambda^1$ resp. the ending point of $\lambda^{q+1}$. We define the collections
    \begin{equation}
    \scrU^-\equiv\{\scrU^-_0,\ldots,\scrU^-_q,\scrV^-\}\;\;\textrm{and}\;\;\scrU^+\equiv\{\scrU^+_1,\ldots,\scrU^+_{q+1},\scrV^+\}.
    \end{equation}
Let $r\in\bbZ_{\geq0}$ with $r+1\leq q+1$. 
\begin{enumerate}
\item We say that 
    \begin{equation}
    \Theta=(\Theta^1,\ldots,\Theta^{r+1})\in\scrW(\Lambda,\scrU^-,\scrU^+)\subset\overline{\scrL}\big(\Pi(\Sigma);f_i-f_j\big)
    \end{equation}
if there exists non-negative integers
    \begin{equation}
    k_0=0<k_1<\cdots<k_r<k_{r+1}=q+1
    \end{equation}
such that the following conditions hold.
\begin{itemize}
\item $\Theta^s\in\scrL(c_{k_{s-1}},c_{k_s};f_i-f_j)$, where $2\leq s\leq r$.
\item $\Theta^s$ exits $\Omega(c_{k_{s-1}})$ through $\scrU^-_{k_{s-1}}$ and enters $\Omega(c_{k_s})$ through $\scrU^+_{k_s}$, where $2\leq s\leq r$.
\item $\Theta^1\in\scrL\big(\Pi(\Sigma),c_1;f_i-f_j\big)$.
\item $\Theta^1$ starts at a point in $\scrU^-_0$ and enters $\Omega(c_{k_1})$ through $\scrU^+_{k_1}$.
\item $\Theta^{r+1}\in\scrL\big(c_{k_r},\Pi(\Sigma);f_i-f_j\big)$.
\item $\Theta^{r+1}$ exits $\Omega(c_{k_r})$ through $\scrU^-_{k_r}$ and ends at a point in $\scrU^+_{q+1}$.
\end{itemize}
\item We say that 
    \begin{equation}
    \theta\equiv(\theta^1,\ldots,\theta^{r+1})\in\scrW(\lambda,\scrU^-,\scrU^+)\subset\overline{\scrL}_0\big(\Pi(\Sigma);f_i-f_j\big)
    \end{equation}
if the following conditions hold, where $\Theta$ is the maximal extension of $\theta$.
\begin{itemize}
\item $\Theta\in\scrW(\Lambda,\scrU^-,\scrU^+)$.
\item $\theta^1$'s starting point is in $\scrV^-$ and $\theta^{r+1}$'s ending point is in $\scrV^+$.
\end{itemize}
\end{enumerate}
The collection 
    \begin{equation}
    \big\{\scrW(\lambda,\scrU^-,\scrU^+)\big\}\;\;\textrm{resp.}\;\;\big\{\scrW(\Lambda,\scrU^-,\scrU^+)\big\}
    \end{equation}
provides a fundamental system of open neighborhoods for a metrizable topology on 
    \begin{equation}
    \overline{\scrL}_0\big(\Pi(\Sigma);f_i-f_j\big)\;\;\textrm{resp.}\;\;\overline{\scrL}\big(\Pi(\Sigma);f_i-f_j\big).
    \end{equation}

Finally, we define both $\overline{\scrL}_0(f_i-f_j)$ and $\overline{\scrL}(f_i-f_j)$ to have the topology induced by the disjoint union.

\subsubsection{Compactness}
In this part, we prove the following result.

\begin{prop}\label{prop:edgeconvergence}
$\overline{\scrL}(f_i-f_j)$ is compact.
\end{prop}

\begin{proof}
It suffices to prove the following statement: for any sequence 
    \begin{equation}
    \{\Lambda_\nu\}\subset\overline{\scrL}(f_i-f_j),
    \end{equation}
there exists a subsequence which converges. Immediately, we may reduce to the case that each $\Lambda_\nu$ is actually an unbroken maximally extended gradient flow $\Phi_\nu$ of $f_i-f_j$; the general case follows by considering each associated sequence of components separately. Moreover, by passing to a subsequence, we may assume each $\Phi_\nu$ is of the same type, i.e., each $\Phi_\nu$ is either: 
\begin{enumerate}
\item a maximally extended Morse flow connecting $a=c_0$ to $b=c_{q+1}$,
\item a maximally extended fold terminating flow starting at $a=c_0$, 
\item a maximally extended fold emanating flow ending at $b=c_{q+1}$, 
\item or a maximally extended singular flow.
\end{enumerate}
For clarity, we find it best to handle each case separately, even though the differences in the arguments are minimal.
\begin{enumerate}
\item Let $p^-_{\nu,0}\in\partial\Omega(c_0)$ be the exit point of $\Phi_\nu$ from a Morse neighborhood $\Omega(c_0)$ of $c_0$. Since 
    \begin{equation}
    W^u(c_0)\cap\partial\Omega(c_0)\cong S^{I(c_0)-1},
    \end{equation}
we may pass to a subsequence and assume $\{p^-_{\nu,0}\}$ converges to a point 
    \begin{equation}
    p^-_0\in W^u(c_0)\cap\partial\Omega(c_0);
    \end{equation}
let $\Lambda^1$ be the unique maximally extended gradient flow of $f_i-f_j$ at $p^-_0$.

We claim there exists $c_1\in\crit(f_i-f_j)$ such that $\Lambda^1\in\calL(c_0,c_1;f_i-f_j)$. If this were not the case, then $\Lambda^1\in\calL\big(c_0,\Pi(\Sigma);f_i-f_j)$. By dependence of solutions to ODEs on initial conditions, there exists a sequence of points $\{x_\nu\}$ associated to $\{\Phi_\nu\}$ which converges to the ending point of $\Lambda^1$. Using Lemma \ref{lem:prelim2}, we get a contradiction to the assumption that $\{\Phi_\nu\}$ is a sequence of Morse flows.

Let $\Omega(c_1)$ be a Morse neighborhood of $c_1$; we denote by 
    \begin{equation}
    p^+_1\in W^s(c_1)\cap\partial\Omega(c_1)
    \end{equation}
the entry point of $\Lambda^1$ into $\Omega(c_1)$. By dependence of solutions to ODEs on initial conditions, there exists a sequence of entry points
    \begin{equation}
    \{p^+_{\nu,1}\}\subset\partial\Omega(c_1)
    \end{equation}
associated to $\{\Phi_\nu\}$; moreover, by Lemma \ref{lem:prelim1}, $\{p^+_{\nu,1}\}$ converges to $p^+_1$. Now, if $c_1=c_{q+1}$, we are done.

We will refer to the following argument as (*). Suppose $c_1\neq c_{q+1}$, then each $p^+_{1,\nu}\notin W^s(c_1)$. We denote by 
    \begin{equation}
    \{p^-_{\nu,1}\}\in\partial\Omega(c_1)
    \end{equation}
the sequence of exit points associated to $\{\Phi_\nu\}$. Since $\partial\Omega(c_1)$ is compact, we may pass to a subsequence and assume $\{p^-_{\nu,1}\}$ converges to a point 
    \begin{equation}
    p^-_1\in\partial\Omega(c_1).
    \end{equation}
We claim $p^-_1\in W^u(c_1)$. If this were not the case, then there exists a maximally extended gradient flow $\Gamma$ of $f_i-f_j$ that enters $\Omega(c_1)$ through a point $\Gamma^+\in\partial\Omega(c_1)$, satisfying 
    \begin{equation}
    (f_i-f_j)(\Gamma^+)=(f_i-f_j)(p^+_{\nu,1}),\;\;\forall\nu,
    \end{equation}
and exits $\Omega(c_1)$ through a point $\Gamma^-\in\partial\Omega(c_1)$ which is the limit of $\{p^-_{\nu,1}\}$. By Lemma \ref{lem:prelim1}, $\{p^+_{\nu,1}\}$ converges to $\Gamma^+$, i.e., $\Gamma^+=p^+_1\in W^s(c_1)$; but this is an immediate contradiction since $\Gamma^+$ cannot possibly be in $W^s(c_1)$. This concludes the argument labeled (*).

The proof of this case is now completed via induction since (1) $\crit(f_i-f_j)$ is a finite set and (2) $f_i-f_j$ strictly decreases along any gradient flow.
\item Let $p^-_{\nu,0}\in\partial\Omega(c_0)$ be the exit point of $\Phi_\nu$ from a Morse neighborhood $\Omega(c_0)$ of $c_0$. Since 
    \begin{equation}
    W^u(c_0)\cap\partial\Omega(c_0)\cong S^{I(c_0)-1},
    \end{equation}
we may pass to a subsequence and assume $\{p^-_{\nu,0}\}$ converges to a point 
    \begin{equation}
    p^-_0\in W^u(c_0)\cap\partial\Omega(c_0);
    \end{equation}
let $\Lambda^1$ be the unique maximally extended gradient flow of $f_i-f_j$ at $p^-_0$.

If $\Lambda^1\in\calL\big(c_0,\Pi(\Sigma);f_i-f_j)$, then we are done by Lemma \ref{lem:prelim2}. Otherwise, there exists $c_1\in\crit(f_i-f_j)$ such that $\Lambda^1\in\calL(c_0,c_1;f_i-f_j)$ -- we now perform the argument (*).

The proof of this case is now completed via induction since (1) $\crit(f_i-f_j)$ is a finite set and (2) $f_i-f_j$ strictly decreases along any gradient flow.
\item Let $p^-_{\nu,0}\in\partial U$ be the starting point of $\Phi_\nu$; by passing to a subsequence, we may assume $\{p^-_{\nu,0}\}$ converges to a point $p^-_0\in\partial U$. By the local from of $\grad(f_i-f_j)$ near the cusp-edge, we see that there exists a unique maximally extended gradient flow $\Lambda^1$ starting at $p^-_0$.

We claim there exists $c_1\in\crit(f_i-f_j)$ such that $\Lambda^1\in\calL\big(\Pi(\Sigma),c_1;f_i-f_j\big)$. If this were not the case, then $\Lambda^1\in\calL\big(\Pi(\Sigma);f_i-f_j\big)$. By dependence of solutions to ODEs on initial conditions, there exists a sequence of points $\{x_\nu\}$ associated to $\{\Phi_\nu\}$ which converges to the ending point of $\Lambda^1$. Using Lemma \ref{lem:prelim2}, we get a contradiction to the assumption that $\{\Phi_\nu\}$ is a sequence of fold terminating flows.

We may now perform the argument (*). The proof of this case is now completed via induction since (1) $\crit(f_i-f_j)$ is a finite set and (2) $f_i-f_j$ strictly decreases along any gradient flow.
\item Let $p^-_{\nu,0}\in\partial U$ be the starting point of $\Phi_\nu$; by passing to a subsequence, we may assume $\{p^-_{\nu,0}\}$ converges to a point $p^-_0\in\partial U$. By the local from of $\grad(f_i-f_j)$ near the cusp-edge, we see that there exists a unique maximally extended gradient flow $\Lambda^1$ starting at $p^-_0$.

If $\Lambda^1\in\calL\big(\Pi(\Sigma);f_i-f_j)$, then we are done by Lemma \ref{lem:prelim2}. Otherwise, there exists $c_1\in\crit(f_i-f_j)$ such that $\Lambda^1\in\calL\big(\Pi(\Sigma),c_1;f_i-f_j\big)$ -- we may now perform the argument (*).

The proof of this case is now completed via induction since (1) $\crit(f_i-f_j)$ is a finite set and (2) $f_i-f_j$ strictly decreases along any gradient flow.
\end{enumerate}
This completes the proof of all cases; the proposition follows.
\end{proof}

An immediately corollary of the previous proposition is the following.

\begin{cor}\label{cor:edgeconvergence}
$\overline{\scrL}_0(f_i-f_j)$ is compact.
\end{cor}

\section{Compactifying stable Morse flow trees}
\subsection{Stable broken Morse flow trees}

Note, given any broken gradient flow $\lambda$ associated to $f_i-f_j$, we may define its 1-jet lift 
    \begin{equation}
    \Big(\widetilde{\lambda}^1,\widetilde{\lambda}^2\Big)
    \end{equation}
in the natural way by lifting each component. Moreover, this will be an ordered pair since we require that $\widetilde{\lambda}^1$ resp. $\widetilde{\lambda}^2$ lies in the sheet determined by $f_i$ resp. $f_j$. We may analogously define its cotangent lift
    \begin{equation}
    \Big(\overline{\lambda}^1,\overline{\lambda}^2\Big).
    \end{equation}

\begin{defin}\label{defin:prestablemorseflowtree}
A \emph{prestable Morse flow tree} of $L$ is a source tree $T$ decorated with the following data.
\begin{enumerate}
\item For each $e\in\edgef(T)$, we associate a non-constant broken gradient flow $\lambda$ of $L$.
\item Let $v\in\vertex(T)$ be a $k$-valent vertex with cyclically ordered edges $e_1,\ldots,e_k$. Recall, associated to $e_j$ is its cotangent lift 
    \begin{equation}
    \Big(\overline{\lambda}^1_j,\overline{\lambda}^2_j\Big).
    \end{equation}
We require the following equality: 
    \begin{equation}
    \overline{p}_j\equiv\overline{\lambda}^{in}_j(v)=\overline{\lambda}^{out}_{j+1}(v)\in\Pi_\bbC(L).
    \end{equation}
\item We require that all of the cotangent lifts of all of the edges of $T$ fit together to give an oriented loop in $\Pi_\bbC(L)$.
\end{enumerate}
\end{defin}

\begin{rem}
Observe, we have the following kinds of 2-valent vertices allowed in property (2) of Definition \ref{defin:prestablemorseflowtree}. It may be the case that $v$ is a 2-valent vertex such that the broken gradient flow $\lambda_1$ of $L$ associated to $e_1$ and the broken gradient flow $\lambda_2$ of $L$ associated to $e_2$ are both broken gradient flows of $L$ associated to $f_i-f_j$ and $v\notin\crit(f_i-f_j)$. We call these kinds of 2-valent vertices \emph{removable}. Removable vertices are necessary for (pre)-Floer-Gromov convergence as they can be approached by collapsing edges in a sequence of (pre)stable Morse flow trees, cf. Figure \ref{fig:ekholm19} where a removable vertex is approached by the collision of a $Y_0$-vertex with an adjacent switch and end.
\end{rem}

\begin{defin}\label{defin:stablemorseflowtree}
A \emph{stable Morse flow tree} of $L$ is a source tree $T$ decorated with the following data.
\begin{enumerate}
\item For each $e\in\edgef(T)$, we associate a non-constant gradient flow $\phi$ of $L$.
\item Let $v\in\vertex(T)$ be a $k$-valent vertex with cyclically ordered edges $e_1,\ldots,e_k$. Recall, associated to $e_j$ is its cotangent lift 
    \begin{equation}
    \Big(\overline{\phi}^1_j,\overline{\phi}^2_j\Big).
    \end{equation}
We require the following equality: 
    \begin{equation}
    \overline{p}_j\equiv\overline{\phi}^{in}_j(v)=\overline{\phi}^{out}_{j+1}(v)\in\Pi_\bbC(L).
    \end{equation}
Moreover, we require that $v$ is not a removable vertex.
\item We require that all of the cotangent lifts of all of the edges of $T$ fit together to give an oriented loop in $\Pi_\bbC(L)$.
\end{enumerate}
\end{defin}

The action of \emph{stabilizing} a prestable Morse flow trees means adding 2-valent vertices along critical points of broken gradient flows and removing removable vertices. The action of \emph{destabilizng} a (pre)stable Morse flow tree means removing 2-valent vertices along critical points of broken gradient flows and adding removable vertices. Clearly, any prestable Morse flow tree can be stabilized to a stable Morse flow tree and any stable Morse flow tree can be destablized to a prestable Morse flow tree. However, for any stable Morse flow tree, there is not a unique prestable Morse flow tree which stabilizes to it.

As with Morse flow trees, (pre)stable Morse flow trees admit nodes.

\begin{defin}\label{defin:node2}
Let $v\in\vertex(T)$ be a vertex of a (pre)stable Morse flow tree $T$ and consider a non-trivial decomposition
     \begin{equation}
    I\amalg J=\{e\in\edgef(T):v\in e\},
    \end{equation}
where both $I$ and $J$ consist of consecutively ordered edges. We say $v$ has a \emph{node} if there exists a decomposition 
    \begin{equation}
    T_I\cup_{v_I\sim v_J} T_J=T,
    \end{equation}
where $I$ and $J$ are as above, such that $T_I$ resp. $T_J$ is a (pre)stable Morse flow tree containing a vertex $v_I$ resp. $v_J$ satisfying
    \begin{equation}
    I=\{e\in\edgef(T_I):v_I\in e\}\;\;\mathrm{resp.}\;\;J=\{e\in\edgef(T_J):v_J\in e\}.
    \end{equation}
We say $T$ is \emph{non-nodal} if it contains no nodes.
\end{defin}

\begin{lem}
There exists a decomposition of $T$ into non-nodal (pre)stable Morse flow trees:
    \begin{equation}
    T=T_{I_1}\cup_{v_{I_1}\sim v_{I_2}'}T_{I_2}\cup_{v_{I_2}\sim v_{I_2}'}\cup\cdots\cup_{v_{I_{s-1}}\sim v_{I_s}'} T_{I_s}.
    \end{equation}
Moreover, this decomposition is unique up to relabeling/reordering.
\end{lem}

\begin{rem}\label{rem:node2}
Observe, any 2-valent vertex $v$ in a stable Morse flow tree such that the gradient flow $\phi_1$ of $L$ associated to $e_1$ and the gradient flow $\phi_2$ of $L$ associated to $e_2$ are both broken gradient flows of $L$ associated to $f_i-f_j$ and $v\in\crit(f_i-f_j)$ is a node. Of course, there are other kinds of nodes as in Remark \ref{rem:node1}.
\end{rem}

Let $T$ be a (pre)stable Morse flow tree; we may associate a (non-unique) stable metric ribbon tree $T^\cs$ to $T$ by (1) adding semi-infinite edges to various vertices, (2) choosing an arbitrary root from among the aforementioned added semi-infinite edges, and (3) choosing a metric such that the length of an interior edge is the length of the associated (broken) gradient flow of $L$ (with respect to the metric $g$ on $M$). We will refer to $T^\cs$ as a \emph{combinatorial stabilization} of $T$. Clearly, a combinatorial stabilization of a (pre)stable Morse flow tree is not unique; however, all of the discussion below is evidently independent of this auxiliary choice.

The following is the main technical definition of the present article.

\begin{defin}\label{defin:prefloergromov}
Let $\{T_\nu\}$ be a sequence of prestable Morse flow trees and $T$ a prestable Morse flow tree. We say $\{T_\nu\}$ \emph{pre-Floer-Gromov converges} to $T$ if the following conditions hold.
\begin{enumerate}
\item There exists a sequence of combinatorial stabilizations $\{T^\cs_\nu\}$ associated to $\{T_\nu\}$ and a combinatorial stabilization $T^\cs$ of $T$ such that $\{T^\cs_\nu\}$ converges to $T^\cs$ in the moduli space of stable metric ribbon trees. This condition already has the following immediate consequences. 
    \begin{enumerate}
    \item For $\nu$ sufficiently large, the domain of each $T_\nu$ is isomorphic as a finite tree to a fixed finite tree $\Gamma$.
    \item $T^\cs$ is a stable metric ribbon tree such that (1) the underlying finite tree of $T^\cs$ is isomorphic as a finite tree to $\Gamma$ and (2) the metric is such that, after collapsing each edge of the underlying finite tree of $T^\cs$ of length 0, we obtain a finite tree isomorphic as a finite tree to the domain of $T$.
    \end{enumerate}
\item For $\nu$ sufficiently large, we may consider a sequence $\{e_\nu\}$ of interior edges of $\{T_\nu\}$ corresponding to the same edge in $\Gamma$ and the corresponding interior edge $e$ of $T^\cs$. We require that each broken gradient flow $\lambda_\nu$ of $L$ associated to each $e_\nu$ are all actually broken gradient flows of $L$ in the same moduli space $\overline{\calL}_0(f_i-f_j)$. We also require that the broken gradient flow $\lambda$ of $L$ associated to $e$ is in $\overline{\calL}_0(f_i-f_j)$ (if $e$ is length 0, then $\lambda$ is constant). Finally, we require that $\{\lambda_\nu\}$ converges to $\lambda$. Moreover, we require the aforementioned conditions for all such sequence of corresponding interior edges.
\end{enumerate}
\end{defin}

\begin{defin}\label{defin:floergromov}
Let $\{T_\nu\}$ be a sequence of stable Morse flow trees and $T$ a stable Morse flow tree. We say $\{T_\nu\}$ \emph{Floer-Gromov converges} to $T$ if there exists a sequence of prestable Morse flow trees $\{T^\pre_\nu\}$ which stabilizes to $\{T_\nu\}$ and a prestable Morse flow tree $T^\pre$ which stabilizes to $T$ such that $\{T^\pre_\nu\}$ pre-Floer-Gromov converges to $T^\pre$.
\end{defin}

In the sequel we will require the following proposition whose proof we collect here. 

\begin{prop}\label{prop:convergence}
Let $\{T_\nu\}$ be a sequence of stable Morse flow trees together with a constant $N\in\bbZ_{>0}$ such that the domain of any $T_\nu$ has at most $N$ edges. Perhaps after passing to a subsequence, there exists a stable Morse flow tree $T$ such that $\{T_\nu\}$ Floer-Gromov converges to $T$.
\end{prop}

\begin{proof}
We consider the associated sequence $\{T^\pre_\nu\}$. First, since 
    \begin{equation}
    \sup_{e\in\edgef(T_\nu)}T_\nu=N<\infty,
    \end{equation}
we may pass to a subsequence and assume each domain of $T^\pre_\nu$ is isomorphic as a finite tree to a fixed finite tree $\Gamma$. Moreover, by passing to a subsequence we may assume there exists a constant $k_N\in\bbZ_{>0}$ and a sequence $\{T^{\pre,\cs}_\nu\}$ of combinatorial stabilizations associated to $\{T^\pre_\nu\}$ such that 
    \begin{equation}
    \{T^{\pre,\cs}_\nu\}\subset\overline{G}_{k_N}
    \end{equation}
converges to some $T^{\pre,\cs}\in\overline{G}_{k_N}$. Observe, the underlying finite tree of $T^{\pre,\cs}$ is isomorphic as a finite tree to $\Gamma$. Consider any associated sequence of corresponding edges $\{e_\nu\}$ of $\{T^\pre_\nu\}$. Since there are only finitely many local function differences, we may pass to a subsequence and assume the sequence $\{\lambda_\nu\}$ of broken gradient flows of $L$ associated to $\{e_\nu\}$ satisfies 
    \begin{equation}
    \{\lambda_\nu\}\subset\overline{\calL}_0(f_i-f_j)
    \end{equation}
for some choice of $f_i-f_j$. Thus, by Corollary \ref{cor:edgeconvergence}, we may pass to a subsequence and assume $\{\lambda_\nu\}$ converges to a broken gradient flow $\lambda$ of $f_i-f_j$. We assign this broken gradient flow to the edge $e$ of $\Gamma$ associated to the sequence $\{e_\nu\}$. Since, again, there is a uniform bound on the number of edges of the domains of $\{T_\nu\}$, we may pass to a subsequence and assume the previous claims are true for each associated sequence of corresponding edges of $\{T^\pre_\nu\}$. Note, if $\lambda$ has length 0, then $e$ has length 0 when thought of as an interior edge of $T^{\pre,\cs}$. In particular, by collapsing all edges of the underlying finite tree of $T^{\pre,\cs}$ of length 0, we obtain a prestable Morse flow tree $T^\pre$ whose stabilization $T$ is the Floer-Gromov limit of $\{T_\nu\}$ by construction.
\end{proof}

\begin{rem}\label{rem:boundededges}
We still consider the setup of Proposition \ref{prop:convergence}. It should be readily apparent that, since: (1) there are finitely many local function differences and (2) each local function difference has finitely critical points; there exists a constant $k(M,L,N)\in\bbZ_{>0}$, depending on $M$, $L$, and $N$, such that the domain of $T$ has at most $k(M,L,N)$ edges.
\end{rem}

\subsection{Floer-Gromov topology}
Let $\frakS(L)$ be the moduli space of stable Morse flow trees of $L$ and $\frakS(L,N)$ the moduli space of stable Morse flow trees of $L$ whose domains have at most $N\in\bbZ_{>0}$ edges. In this subsection, we will construct the Floer-Gromov topology on the moduli space 
    \begin{equation}
    \overline{\frakS}(L,N)\equiv\big\{T\in\frakS(L):\exists\{T_\nu\}\subset\frakS(L,N)\;\;\textrm{Floer-Gromov converging to $T$}\big\}
    \end{equation}
and prove Theorem \ref{thm:main}. 

First, we will take a slight digression to explain some preliminaries on convergence structures, cf. \cite[Section 5.6]{MS12}.

\begin{defin}
Let $X$ be a set. A \emph{convergence structure} on $X$ is a subset 
    \begin{equation}
    \calC\subset X\times X^\bbN
    \end{equation}
satisfying the following axioms. (We use the notation $(x_\nu)_\nu$ for the map $\nu\mapsto x_\nu$.)
\begin{itemize}
\item\emph{Constant:} If $\big(x,(x_\nu)_\nu\big)$ satisfies $x_\nu=x$ for all $\nu$, then $\big(x,(x_\nu)_\nu\big)\in\calC$.
\item\emph{Subsequence:} If $\big(x,(x_\nu)_\nu\big)\in\calC$ and $g:\bbN\to\bbN$ is strictly increasing, then $\big(x,(x_{g(\nu)})_\nu\big)\in\calC$.
\item\emph{Subsubsequence:} If for every strictly increasing function $g:\bbN\to\bbN$ there exists a strictly increasing function $f:\bbN\to\bbN$ such that $\big(x,(x_\nu)_{(g\circ f)(\nu)}\big)\in\calC$, then $\big(x,(x_\nu)_\nu\big)\in\calC$.
\item\emph{Uniqueness:} If $\big(x,(x_\nu)_\nu\big),\big(y,(x_\nu)_\nu\big)\in\calC$, then $x=y$.
\item\emph{Diagonal:} If $\big(x,(x_\kappa)_\kappa\big)\in\calC$ and, for every $\kappa$, $\big(x_\kappa,(x_{\nu,\kappa})_\nu\big)\in\calC$, then there exists subsequences $\{\nu_r\}_{r\in\bbN}$ and $\{\kappa_r\}_{r\in\bbN}$ such that $\big(x,(x_{\nu_r,\kappa_r})_r\big)\in\calC$.
\end{itemize}
\end{defin}

Intuitively, we should think of an element $\big(x,(x_\nu)_\nu\big)$ as a sequence $\{x_\nu\}$ that converges to $x$; this is made precise by the following definition and lemma.

\begin{defin}\label{defin:csopen}
Let $(X,\calC)$ be a convergence space. We define a topology on $X$ as follows; $U\subset X$ is open if and only if, for every 
    \begin{equation}
    \big(x,(x_\nu)_\nu\big)\in\calC\cap\Big(U\times X^\bbN\Big),
    \end{equation}
it follows that $x_\nu\in U$ for $\nu$ sufficiently large.
\end{defin}

\begin{lem}
The collection of open sets in Definition \ref{defin:csopen} actually forms a topology on $X$. Moreover, this topology is unique with respect to the property that a sequence $\{x_\nu\}\subset X$ converges to $x\in X$ if and only if $\big(x,(x_\nu)_\nu\big)\in\calC$.
\end{lem}

Returning to stable Morse flow trees, we consider the subset   
    \begin{equation}
    \calC_\mathrm{FG}\subset\overline{\frakS}(L,N)\times\overline{\frakS}(L,N)^\bbN
    \end{equation}
consisting of the pairs 
    \begin{equation}
    \Big\{\big(T,(T_\nu)_\nu\big):\{T_\nu\}\;\;\textrm{Floer-Gromov converges to $T$}\Big\}.
    \end{equation}

\begin{lem}\label{lem:topology}
$\big(\overline{\frakS}(L,N),\calC_\mathrm{FG}\big)$ is a convergence space.
\end{lem}

\begin{proof}
Of course, we must verify the \emph{constant}, \emph{subsequence}, \emph{subsubsequence}, \emph{uniqueness}, and \emph{diagonal} axioms. These are all straightforward since Floer-Gromov convergence reduces to convergence of (1) stable metric ribbon trees and (2) (broken) gradient flows.\footnote{Moreover, the uniqueness axiom will also require the fact that we consider stable Morse flow trees -- uniqueness would not hold for prestable Morse flow trees essentially due to removable vertices and the fact that broken gradient flows may decorate edges.} We will prove the \emph{diagonal} axiom to illustrate the argument. 

Let $\{T_{\nu}\}\subset\overline{\frakS}(L,N)$ be a sequence Floer-Gromov converging to $T\in\overline{\frakS}(L,N)$. Moreover, for each fixed $\nu$, let $\{T_{\nu,\kappa}\}\subset\overline{\frakS}(L,N)$ be a sequence Floer-Gromov converging to $T_\nu$. We work with the associated $T^\pre$, $\{T^\pre_\nu\}$, and $\{T^\pre_{\nu,\kappa}\}$. (The follow discussion is technically only true for $\nu,\kappa\gg0$, but we will elide this point for brevity.)
\begin{enumerate}
\item By definition, there exists an associated sequence of combinatorial stabilizations $\{T^{\pre,\cs}_\nu\}$ and a combinatorial stabilization $T^{\pre,\cs}$ such that $\{T^{\pre,\cs}_\nu\}$ converges to $T^{\pre,\cs}$ in the moduli space of stable metric ribbon trees. Analogously, for each fixed $\nu$, we define $\{T^{\pre,\cs}_{\nu,\kappa}\}$ and $\underline{T}^{\pre,\cs}_\nu$. The main observation to make is that the definition of pre-Floer-Gromov convergence implies
    \begin{equation}\label{eqn:equality}
    \underline{T}^{\pre,\cs}_\nu=T^{\pre,\cs}_\nu,\;\;\forall\nu
    \end{equation}
in the moduli space of stable metric ribbon trees; in particular, after collapsing a subset of length 0 interior edges of $\underline{T}^{\pre,\cs}_\nu$, we see the aforementioned equality. This immediately implies that there exists a subsequence $\{(\nu_r,\kappa_r)\}$ such that $\{T^{\pre,\cs}_{\nu_r,\kappa_r}\}$ converges to $T^{\pre,\cs}$.
\item We now make a similar argument for the associated sequences of corresponding edges. I.e., for any fixed $\nu$, consider any sequence $\{e_{\nu,\kappa}\}$ of corresponding interior edges associated to $\{T^{\pre,\cs}_{\nu,\kappa}\}$ and the corresponding interior edge $\underline{e}_\nu$ of $\underline{T}^{\pre,\cs}_\nu$ such that the sequence of broken gradient flows $\{\lambda_{\nu,\kappa}\}$ associated to $\{e_{\nu,\kappa}\}$ converges to the broken gradient flow $\underline{\lambda}_\nu$ associated to $e_\nu$. There are two cases.
    \begin{enumerate}
    \item First, for each $\nu$, $\underline{e}_\nu$ and $\underline{\lambda}_\nu$ are length 0. Hence, the sequence $\{\underline{\lambda}_\nu\}$ consists of constant broken gradient flows and, perhaps after passing to a subsequence, converges to the image of a vertex of $T$.
    \item Second, there exists a sequence $\{e_\nu\}$ of corresponding interior edges associated to $\{T^{\pre,\cs}_\nu\}$ and a corresponding interior edge $e$ of $T^{\pre,\cs}$, where the sequence of broken gradient flows $\{\lambda_\nu\}$ associated to $\{e_\nu\}$ converges to the broken gradient flow $\lambda$ associated to $e$, such that, using the equality \eqref{eqn:equality}, $\underline{e}_\nu=e_\nu$ and $\underline{\lambda}_\nu=\lambda_\nu$. This immediately implies there exists a subsequence $\{(\nu_r,\kappa_r)\}$ such that $\{\lambda_{\nu_r,\kappa_r}\}$ converges to $\lambda$. 
    \end{enumerate}
We now induct over the number of edges.
\end{enumerate}
Hence, there exists a subsequence $\{(\nu_r,\kappa_r)\}$ such that $\{T_{\nu_r,\kappa_r}\}$ Floer-Gromov converges to $T$.
\end{proof}

Well call the topology on $\overline{\frakS}(L,N)$ induced by the convergence structure $\calC_\mathrm{FG}$ the \emph{Floer-Gromov topology}. We are now ready to prove our main result.

\begin{proof}[Proof of Theorem \ref{thm:main}]
Of course, it suffices to prove the following statement: for any sequence 
    \begin{equation}
    \{T_\nu\}\subset\overline{\frakS}(L,N), 
    \end{equation}
there exists a subsequence which Floer-Gromov converges. There are really two parts of the proof. 
\begin{enumerate}
\item First, showing that, perhaps after passing to a subsequence, there exists an element $T\in\frakS(L)$ such that $\{T_\nu\}$ Floer-Gromov converges to $T$.

\item Second, showing that $T$ is actually an element of $\overline{\frakS}(L,N)$.
\end{enumerate}

For the first part we observe that, for each fixed $\nu$, there exists a sequence $\{T_{\nu,\kappa}\}\subset\frakS(L,N)$ such that $\{T_{\nu,\kappa}\}$ Floer-Gromov converges to $T_\nu$ by definition. In particular, by Remark \ref{rem:boundededges}, the domain of $T_\nu$ has at most $k(M,L,N)$ edges; thus, the first part follows by Proposition \ref{prop:convergence}.

The second part is straightforward. We wish to show there exists a sequence in $\frakS(L,N)$ which converges to the $T$ constructed in the first part; by repeating the argument proving the \emph{diagonal} axiom in Lemma \ref{lem:topology}, we may find a subsequence $\{(\nu_r,\kappa_r)\}$ such that $\{T_{\nu_r,\kappa_r}\}$, which is a subsequence of the sequence from the first part in $\frakS(L,N)$, Floer-Gromov converges to $T$.
\end{proof}

\subsection{An example}
We will end by reviewing Ekholm's examples of (pre)stable Morse flow trees and describe their Floer-Gromov limits, cf. \cite[Section 7]{Ekh07}.

\begin{figure}[h]
\centering
\scalebox{.5}{\includegraphics{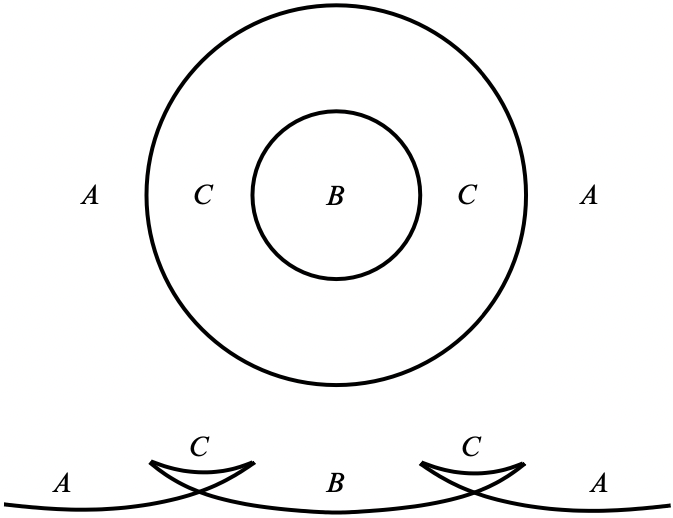}}
\caption{(Figure 17 in \cite{Ekh07})  (Top) The top view of $F_1'$. (Bottom) The cross-sectional view of $F_1'$.}\label{fig:ekholm17}
\end{figure}

Consider two fronts $F_0,F_1\subset J^0\bbR^2\cong\bbR^3$, where $F_0$ is the 0-section and $F_1$ is the graph of the function $z(x_1,x_2)=K-x_1$ with $K\in\bbR$ a sufficiently large constant. We will consider the projection of the fronts to the base $\bbR^2$ in a square
    \begin{equation}
    -R\leq\abs{x_i}\leq R,\;\;R\in\bbR_{>0}.
    \end{equation}
We Legendrian isotope $F_1$ into $F_1'$, where $F_1'$ is a front with three sheets, denoted $A$, $B$, and $C$, whose projection to the base $\bbR^2$ looks as follows; we have two concentric circles in $\bbR^2$ and: $A$ is the outer region, $B$ is the inner region, and $C$ is the annular region. Also, we denote by $D$ the single sheet of $F_0$. See Figure \ref{fig:ekholm17}.

For consistency, we will describe our Morse flow trees using Ekholm's notation. Vertices will be labeled by Greek letters and edges will be labeled by positive numbers. For an edge labeled $m\in\bbZ_{>0}$, we write $m\equiv X\vert Y$, where $X,Y\in\{A,B,C,D\}$, if the edge $m$ has a gradient flow associated to the sheets $X$ and $Y$.

\begin{figure}[h]
\centering
\scalebox{.7}{\includegraphics{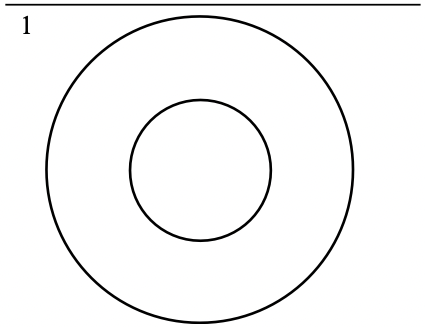}}
\caption{(Figure 18 in \cite{Ekh07}) The sequence of stable Morse flow trees $\{T_\nu\}$.}\label{fig:ekholm18}
\end{figure}

\begin{figure}[h]
\centering
\scalebox{.7}{\includegraphics{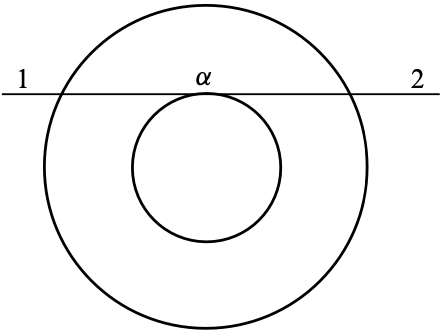}}
\caption{(Figure 19 in \cite{Ekh07}) The prestable Morse flow tree $T^\pre$.}\label{fig:ekholm19}
\end{figure}

We will describe our (pre)stable Morse flow trees using the top view of $F_1'$. First, consider the sequence $\{T_\nu\}$ of one-edge stable Morse flow trees, with edge $1=A\vert D$, starting in the region $A$ and flowing horizontally to the region $A$; $T_1$ flows while remaining in the region $A$ and the sequence moves this edge downwards such that the Floer-Gromov limit is the one-edge stable Morse flow tree $T$, with edge $1=A\vert D$, which is tangent to the circle connecting the region $C$ to the region $B$. We will consider the two-edge prestable Morse flow tree $T^\pre$, with edges $1=2=A\vert D$ and vertex $\alpha$ mapping to the tangency point, which stabilizes to $T$ (note, $\alpha$ is not a puncture). See Figures \ref{fig:ekholm18} and \ref{fig:ekholm19}.

\begin{figure}[h]
\centering
\scalebox{.7}{\includegraphics{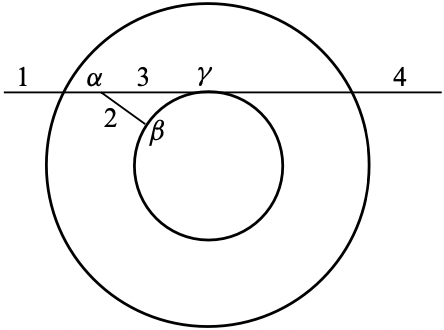}}
\caption{(Figure 20 in \cite{Ekh07}) The sequence of stable Morse flow trees $\{T_\nu'\}$.}\label{fig:ekholm20}
\end{figure}

Now, consider the sequence $\{T_\nu'\}$ of four-edge stable Morse flow trees with $T_1$ described as follows. An edge $1=A\vert D$, flowing from the region $A$ to the region $C$, which ends at a $Y_0$-vertex $\alpha$. From $\alpha$ we have two other edges; first, an edge $2=A\vert C$, flowing transversely into the circle connecting the region $C$ to the region $B$, ending at an end-vertex $\beta$; second, an edge $3=C\vert D$, flowing tangentially into the circle connecting the region $C$ to the region $B$, which ends at a switch-vertex $\gamma$. Finally, from $\gamma$ we have one other edge $4=A\vert D$, flowing from the region $C$ to the region $A$. See Figure \ref{fig:ekholm20}. The sequence shrinks the edges 2 and 3, therefore colliding the vertices $\alpha$, $\beta$, and $\gamma$. In particular, the Floer-Gromov limit of $\{T_\nu'\}$ is our $T$ from above since $\{T_\nu'\}$ pre-Floer-Gromov converges to our $T^\pre$ from above.

\bibliography{References}{}
\bibliographystyle{alpha.bst}
\end{document}